\documentclass[12pt,reqno]{amsart}
\usepackage{latexsym, amsfonts, amsmath, amsthm, amssymb,amscd,color}
\usepackage[dvips]{graphicx}
\usepackage{eepic}
\usepackage{color}
\newtheorem{thm}{Theorem}[section]
\newtheorem{rmk}{Remark}[section]
\newtheorem{prop}[thm]{Proposition}

\newtheorem{cor}[thm]{Corollary}

\newtheorem{problem}[thm]{ Problem}

\numberwithin{equation}{section}
\def\pf{\noindent {\it Proof.} }

\def\1{$\bf{1}$}
\def\0{$\bf{0}$}

\def\cro{\rm{cros}}
\def\ne{\rm{nest}}
\def\nec{\rm{ne}}
\def\sec{\rm{se}}
\def\luc{\rm{luc}}
\def\ruc{\rm{ruc}}
\def\le{\rm{left}}
\def\ri{\rm{right}}

\def\int{\rm{int}}

\def\Min{\rm{Min}}

\def\EC{\rm{EC}}
\def\ER{\rm{ER}}

\def\N{\mathcal{N}}
\def\G{\mathcal{G}}
\def\O{\mathcal{O}}
\def\C{\mathcal{C}}

\def\A{\mathcal{A}}
\def\M{\mathcal{M}}
\def\P{\mathcal{P}}
\def\LP{\mathcal{LP}}

\title{ On the symmetry of ascents and descents over 01-fillings of moon polyominoes }
\author{Anisse Kasraoui}

\begin{document}
\maketitle \centerline{\small Universit\'e de Lyon ; \small
Universit\'e Lyon 1} \centerline{\small Institut Camille Jordan CNRS
UMR 5208} \centerline{\small 43, boulevard du 11 novembre 1918}
\centerline{\small F-69622, Villeurbanne Cedex }
\centerline{\small\texttt{anisse@math.univ-lyon1.fr}}

\begin{abstract}
The purpose of this short paper is to put recent results on the
symmetry of the joint distribution of the numbers of crossings and
nestings of two edges over matchings, set partitions  and linked
partitions, in the larger context of the enumeration of increasing
and decreasing chains of length $2$ in fillings of moon polyominoes.
\end{abstract}


\section{Introduction and Main result}

\subsection{Introduction}  Two natural statistics on simple graphs with
vertex set $[n]$ (or any linear order), where as usual
$[n]:=\{1,2,\cdots,n\}$ for any positive integer $n$,
   are the numbers of \emph{crossings and nestings of two edges},
also called \emph{$2$-crossings} and \emph{$2$-nestings}. Let $G$ be
a \emph{simple graph} (no multiple edges and loops) on $[n]$. A
\emph{$2$-crossing} (resp., \emph{$2$-nesting}) of $G$ is just a
sequence of two arcs $(i_1,j_1),(i_2,j_2)$ of $G$ such that
$i_1<i_2<j_1<j_2$ (resp., $i_1<i_2<j_2<j_1$). If we draw the
vertices of $G$ in increasing order on a line and draw the arcs
above the line, crossings and nestings have the obvious geometrical
meaning (see Figure~1 for an illustration). The numbers of
$2$-crossings and $2$-nestings of $G$ will be denoted by $\cro_2(G)$
and $\ne_2(G)$. For instance, if $G$ is the graph represented below
then $\cro_2(G)=4$ and $\ne_2(G)=6$.

\begin{figure}[h!]\label{fig:graph and filling}
\begin{center}
{\setlength{\unitlength}{1mm}
\begin{picture}(50,30)(0,-10)
\put(-2,0){\line(1,0){54}}
\put(0,0){\circle*{1,3}}\put(0,0){\makebox(0,-6)[c]{\small $1$}}
\put(5,0){\circle*{1,3}}\put(5,0){\makebox(0,-6)[c]{\small $2$}}
\put(10,0){\circle*{1,3}}\put(10,0){\makebox(0,-6)[c]{\small $3$}}
\put(15,0){\circle*{1,3}}\put(15,0){\makebox(0,-6)[c]{\small $4$}}
\put(20,0){\circle*{1,3}}\put(20,0){\makebox(0,-6)[c]{\small $5$}}
\put(25,0){\circle*{1,3}}\put(25,0){\makebox(0,-6)[c]{\small $6$}}
\put(30,0){\circle*{1,3}}\put(30,0){\makebox(0,-6)[c]{\small $7$}}
\put(35,0){\circle*{1,3}}\put(35,0){\makebox(0,-6)[c]{\small $8$}}
\put(40,0){\circle*{1,3}}\put(40,0){\makebox(0,-6)[c]{\small $9$}}
\put(45,0){\circle*{1,3}}\put(45,0){\makebox(0,-6)[c]{\small $10$}}
\put(50,0){\circle*{1,3}}\put(50,0){\makebox(0,-6)[c]{\small $11$}}
\qbezier(0,0)(25,16)(40,0)
\qbezier(5,0)(10,5)(10,0)\qbezier(20,0)(22,7)(25,0)
\qbezier(10,0)(20,10)(30,0)\qbezier(25,0)(30,15)(50,0)
\qbezier(25,0)(32.5,7)(40,0)\qbezier(40,0)(42,6)(45,0)
\qbezier(5,0)(10,8)(15,0)
\end{picture}}\hspace{2.5cm}
{\setlength{\unitlength}{0.7mm}
\begin{picture}(50,50)(0,0)
\put(0,0){\line(1,0){55}}\put(0,5){\line(1,0){50}}\put(0,10){\line(1,0){45}}
\put(0,15){\line(1,0){40}}\put(0,20){\line(1,0){35}}\put(0,25){\line(1,0){30}}
\put(0,30){\line(1,0){25}}\put(0,35){\line(1,0){20}}\put(0,40){\line(1,0){15}}
\put(0,45){\line(1,0){10}}\put(0,50){\line(1,0){5}}
\put(0,0){\line(0,1){55}}\put(5,0){\line(0,1){50}}\put(10,0){\line(0,1){45}}
\put(15,0){\line(0,1){40}}\put(20,0){\line(0,1){35}}\put(25,0){\line(0,1){30}}
\put(30,0){\line(0,1){25}}\put(35,0){\line(0,1){20}}\put(40,0){\line(0,1){15}}
\put(45,0){\line(0,1){10}}\put(50,0){\line(0,1){5}}
\put(0,10){\makebox(6,4)[c]{\small \1}}
\put(5,35){\makebox(6,4)[c]{\small \1}}
\put(5,40){\makebox(6,4)[c]{\small \1}}
\put(10,20){\makebox(6,4)[c]{\small \1}}
\put(20,25){\makebox(6,4)[c]{\small \1}}
\put(25,0){\makebox(6,4)[c]{\small \1}}
\put(25,10){\makebox(6,4)[c]{\small \1}}
\put(40,5){\makebox(6,4)[c]{\small \1}}
\put(0,0){\makebox(5,-6)[c]{\tiny $1$}}
\put(5,0){\makebox(5,-6)[c]{\tiny $2$}}
\put(10,0){\makebox(5,-6)[c]{\tiny $3$}}
\put(15,0){\makebox(5,-6)[c]{\tiny $4$}}
\put(20,0){\makebox(5,-6)[c]{\tiny $5$}}
\put(25,0){\makebox(5,-6)[c]{\tiny $6$}}
\put(30,0){\makebox(5,-6)[c]{\tiny $7$}}
\put(35,0){\makebox(5,-6)[c]{\tiny $8$}}
\put(40,0){\makebox(5,-6)[c]{\tiny $9$}}
\put(45,0){\makebox(5,-6)[c]{\tiny $10$}}
\put(51,0){\makebox(5,-6)[c]{\tiny $11$}}
\put(0,0){\makebox(-6,5)[c]{\tiny $11$}}
\put(0,5){\makebox(-6,5)[c]{\tiny $10$}}
\put(0,10){\makebox(-5,5)[c]{\tiny $9$}}
\put(0,15){\makebox(-5,5)[c]{\tiny $8$}}
\put(0,20){\makebox(-5,5)[c]{\tiny $7$}}
\put(0,25){\makebox(-5,5)[c]{\tiny $6$}}
\put(0,30){\makebox(-5,5)[c]{\tiny $5$}}
\put(0,35){\makebox(-5,5)[c]{\tiny $4$}}
\put(0,40){\makebox(-5,5)[c]{\tiny $3$}}
\put(0,45){\makebox(-5,5)[c]{\tiny $2$}}
\put(0,50){\makebox(-5,5)[c]{\tiny $1$}}
\end{picture}}
\end{center}
\caption{The graph
$\{(1,9),(2,3),(2,4),(3,7),(5,6),(6,9),(6,11),(9,10)\}$ and the
corresponding filling of ${\Delta}_{10}$.}
\end{figure}
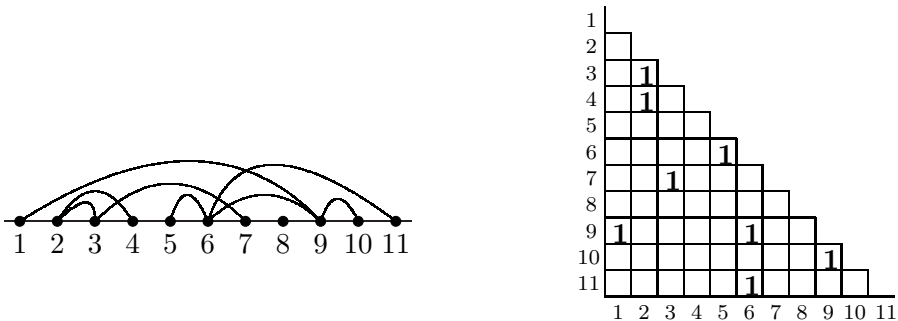

Recently, there has been an increasing interest in studying
crossings and nestings in graphs of matchings, set partitions,
linked partitions and permutations (see e.g.
\cite{Bou,Chen1,Chen2,Co,Mi1,Mi2,KaZe,KlNo,Poz}). In particular, a
property of symmetry has been systematically established. Let
$\M_n$, $\P_n$ and $\LP_n$ be the sets of graphs of matchings, set
partitions and linked partitions of $[n]$, respectively. Then it was
shown that the joint statistic $(\cro_2,\ne_2)$ is symmetrically
distributed over any $\A\in\{\M_n, \P_n, \LP_n\}$ (see respectively
\cite{KlNo,KaZe,Chen1}), i.e. for any integers $k,\ell$,
\begin{equation*}
|\{G\in\A
:\cro_2(G)=k,\ne_2(G)=\ell\}|=|\{G\in\A:\cro_2(G)=\ell,\ne_2(G)=k\}|,
\end{equation*}
or in other words,
\begin{equation}\label{eq:sympartperm}
\sum_{G\in\,\A}p^{\cro_2(G)}q^{\ne_2(G)}=\sum_{G\in\,\A}p^{\ne_2(G)}q^{\cro_2(G)}.
\end{equation}

The purpose of this short paper is to put the latter results in the
more general setting of fillings of arrangements of cells. Note that
such a generalization (for different notions of crossing and
nesting~\cite{Chen2}) was initiated by Krattenthaler~\cite{Kr} and
prolonged by Rubey~\cite{Ru} and De Mier~\cite{Mi1,Mi2}.

 The start is the correspondence between simple
graphs on $[n]$ and $01$-fillings of the triangular Ferrers diagram
$\Delta_n$ (in French notation) of shape $(n-1,n-2,\ldots,1,0)$. The
correspondence consists in labeling columns from top to bottom by
$\{2,3,\ldots,n\}$ and the rows from left to right by
$\{1,2,\ldots,n-1\}$. Then put a 1 in the cell on column labeled~$i$
and  row labeled~$j$ if and only if $(i,j)$ is an arc of $G$. An
illustration is given in Figure~1. In this correspondence, crossings
are send on \emph{descents} and nestings on \emph{ascents}. An
\emph{ascent}, or \emph{North-East chain of length $2$}, (resp.,
\emph{descent} or \emph{South-East chain of length $2$}) in a
$01$-filling is just a set of two 1's in the filling such that one
of them is above and to the right (resp., below and to the left) of
the other and the smallest rectangle contained the two 1's is
contained in the diagram. Then \eqref{eq:sympartperm} can be
translated into a property of symmetry of ascents and descents over
$01$-fillings of the triangular shape $\Delta_n$ according to some
restrictions on the numbers of 1's in columns and rows. In this
paper, we establish such a property in more general arrangements,
namely in moon polyominoes.

\subsection{The main result}

A \emph{polyomino} is an arrangement of square cells. It is
\emph{convex} if along any row of cells and along any column of
cells there is no hole. It is \emph{intersection free} if any two
rows are comparable, i.e., one row can be embedded in the other by
applying a vertical shift. A \emph{moon polyomino} is a convex and
intersection free polyomino. An illustration is given below.
\begin{center}
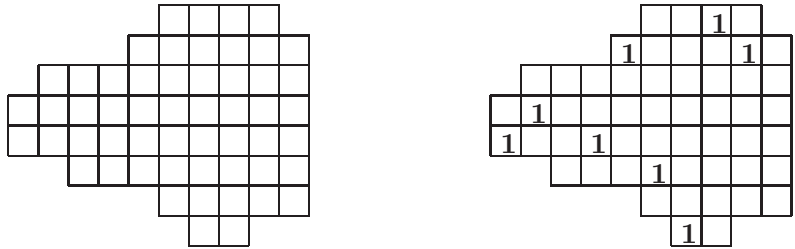
\begin{figure}[h!]
 {\setlength{\unitlength}{0.8mm}
\begin{picture}(50,50)(0,0)
\put(30,0){\line(1,0){10}} \put(25,5){\line(1,0){25}}
\put(10,10){\line(1,0){40}} \put(0,15){\line(1,0){50}}
\put(0,20){\line(1,0){50}} \put(0,25){\line(1,0){50}}
\put(5,30){\line(1,0){45}} \put(20,35){\line(1,0){30}}
\put(15,30){\line(1,0){10}} \put(25,40){\line(1,0){20}}
\put(0,15){\line(0,1){10}}\put(5,15){\line(0,1){15}}
\put(10,10){\line(0,1){20}}\put(15,10){\line(0,1){20}}\put(20,10){\line(0,1){25}}
\put(25,5){\line(0,1){35}}\put(30,0){\line(0,1){40}}\put(35,0){\line(0,1){40}}\put(40,0){\line(0,1){40}}
\put(45,5){\line(0,1){35}} \put(50,5){\line(0,1){30}}
\end{picture}}
\hspace{2cm} {\setlength{\unitlength}{0.8mm}
\begin{picture}(50,50)(0,0)
\put(30,0){\line(1,0){10}} \put(25,5){\line(1,0){25}}
\put(10,10){\line(1,0){40}} \put(0,15){\line(1,0){50}}
\put(0,20){\line(1,0){50}} \put(0,25){\line(1,0){50}}
\put(5,30){\line(1,0){45}} \put(20,35){\line(1,0){30}}
\put(15,30){\line(1,0){10}} \put(25,40){\line(1,0){20}}
\put(0,15){\line(0,1){10}}\put(5,15){\line(0,1){15}}
\put(10,10){\line(0,1){20}}\put(15,10){\line(0,1){20}}\put(20,10){\line(0,1){25}}
\put(25,5){\line(0,1){35}}\put(30,0){\line(0,1){40}}\put(35,0){\line(0,1){40}}\put(40,0){\line(0,1){40}}
\put(45,5){\line(0,1){35}} \put(50,5){\line(0,1){30}}
\put(30,0){\makebox(6,4)[c]{\small \1}}
\put(25,10){\makebox(6,4)[c]{\small \1}}
\put(15,15){\makebox(6,4)[c]{\small\1}}
\put(0,15){\makebox(6,4)[c]{\small \1}}
\put(5,20){\makebox(6,4)[c]{\small\1}}
\put(20,30){\makebox(6,4)[c]{\small\1}}
\put(40,30){\makebox(6,4)[c]{\small \1}}
\put(35,35){\makebox(6,4)[c]{\small \1}}
\end{picture}}
\caption{A moon polyomino $T$ and a filling in $\N(T,{\bf m};A)$,
with ${\bf m}=~(1,2,0,1,2,1,0,1)$ and $A=\{3,10\}$.}
\end{figure}
\end{center}

Let $T$ be a moon polyomino with $s$ rows and $t$ columns. By
convention, we always label the rows from top to bottom by
$\{1,2,\ldots,s\}$ and the columns of $T$ from left to right by
$\{1,2,\ldots,t\}$. A $01$-\emph{filling} $F$ of $T$ consists of
assigning 0 or 1 to each cell. We say that a cell is \emph{empty} if
it has been assigned the value $0$. Given a $s$-uple ${\bf
m}=(m_1,\ldots,m_s)$ (resp., ${\bf m}=(m_1,\ldots,m_t)$) of
integers, we denote by $\N(T,{\bf m})$ (resp., $\N^{\,'}(T,{\bf
m})$) the set of $01$-fillings of $T$ with exactly $m_i$ $1$'s in
row (resp., column) labelled $i$ such that there is at most one $1$
in each column (resp., row). Also, if $A$ is a set of positive
integers, we denote by $\N(T,{\bf m};A)$ (resp., $\N^{\,'}(T,{\bf
m};A)$) the set of fillings $F$ in $\N(T,{\bf m})$ (resp.,
$\N^{\,'}(T,{\bf m}$) whose set of the indices of its empty columns
(resp. rows), denoted $\EC(F)$ (resp., $\ER(F)$), is equal to $A$.
An example is given in Figure~2. Also, if $F$ is a $01$-filling, we
denote by $\nec_2(F)$ and $\sec_2(F)$ the number of ascents and
descents in $F$. For instance, if $F$ is the filling of Figure~1, we
have $\nec_2(F)=6$ and $\sec_2(F)=4$, while for the filling of
Figure~2 we have $\nec_2(F)=\sec_2(F)=4$. We can now state the main
result of the paper, which is a generalization of
\eqref{eq:sympartperm}.

\begin{thm}\label{thm:refinement main1}
For any moon polyomino $T$, the joint statistic $(\nec_2,\sec_2)$ is
symmetrically distributed over each $\N(T,{\bf m};A)$ and each
$\N^{\,'}(T,{\bf m};A)$, i.e., we have
\begin{align*}
\sum_{F\in\,\N(T,{\bf
m};A)}p^{\nec_2(F)}q^{\sec_2(F)}&=\sum_{F\in\,\N(T,{\bf
m};A)}p^{\sec_2(F)}q^{\nec_2(F)},\\
\sum_{F\in\,\N^{\,'}(T,{\bf
m};A)}p^{\nec_2(F)}q^{\sec_2(F)}&=\sum_{F\in\,\N^{\,'}(T,{\bf
m};A)}p^{\sec_2(F)}q^{\nec_2(F)}.
\end{align*}
\end{thm}

Summing over all $A$ in \eqref{eq:symfilling}, we get the following.

\begin{thm}\label{thm: main}
For any moon polyomino $T$, the joint statistic $(\nec_2,\sec_2)$ is
symmetrically distributed over each $\N(T,{\bf m})$ and
$\N^{\,'}(T,{\bf m})$.
\end{thm}

The paper is organized as follows. In Section~2, we give a bijective
proof of Theorem~\ref{thm:refinement main1}. In section~3, we give
an alternative proof of the latter result by computing the
distribution of the joint statistic $(\nec_2,\sec_2)$ over each
$\N(T,{\bf m},A)$. In section~4, we show how we can recover
$\eqref{eq:sympartperm}$ from our result. Finally, we conclude this
paper with some remarks.

 \section{Proof of Theorem~\ref{thm:refinement main1}}

Since the transpose of a moon polyomino is always a moon polyomino,
it suffices to prove the first equation of
Theorem~\ref{thm:refinement main1}. Our first proof is bijective.
Let $T$ be a moon polyomino with $s$~rows and $t$ columns, ${\bf
m}=(m_1,m_2,\ldots,m_s)$ a $s$-uple of nonnegative integers. We will
construct an involution $\Phi:\N(T,{\bf m})\to \N(T,{\bf m})$ such
that for any  $F\in \N(T,{\bf m})$, we have
\begin{align*}
\EC(\Phi(F))=\EC(F),\; \nec_2(\Phi(F))=\sec_2(F),
\;\sec_2(\Phi(F))=\nec_2(F).
\end{align*}

\subsection{Notations}
 The \emph{length-row sequence} of $T$, denoted $r(T)$,
is the sequence $(r_1,r_2,\ldots,r_s)$ where $r_i$ is the length of
$R_i$, the $i$-th row from top. Clearly, the length-row sequence of
any moon polyomino is always unimodal. Thus there exists an unique
$i_0$ such that $r_1\leq r_2\leq\cdots\leq
r_{i_0}>r_{i_0}+1\geq\cdots\geq r_s$. The \emph{upper part of $T$},
denoted $Up(T)$, is the set of rows $R_i$ with $1\leq i\leq i_0$,
and the \emph{lower part}, denoted $Low(T)$, the set of rows $R_i$,
$i_0+1\leq i\leq s$. For instance, if $T$ is the moon polyomino in
Figure~3, we have $r(T)=(4,6,9,10,10,8,5,2)$,
$Up(T)=\{R_1,R_2,R_3,R_4,R_5\}$ and $Low(T)=\{R_6,R_7,R_8\}$. We
order the rows of $T$ by the (total) order $\prec$ defined by
$R_i\prec R_j$ if and only if
\begin{itemize}
\item $r_i<r_j$ or
\item $r_i=r_j$, $R_i\in Up(T)$ and $R_j\in Low(T)$, or
\item $r_i=r_j$, $R_i,R_j\in Up(T)$ and $R_i$ is above $R_j$, or
\item $r_i=r_j$, $R_i,R_j\in Low(T)$ and $R_i$ is below $R_j$.
\end{itemize}
It is clear that $\prec$ is a total order on the rows of $T$.  For
instance, if $T$ is the moon polyomino in Figure~3 we have $R_8\prec
R_1\prec R_7\prec R_2\prec R_6\prec R_3\prec R_4\prec R_5$.

For $i$ an integer, $1\leq i\leq s$, the \emph{$i$-th rectangle of
$T$}, is the greatest rectangle contained in $T$ whose top (resp.,
bottom) row is $R_i$ if $R_i\in Up(T)$ (resp., $R_i\in Low(T)$). An
illustration is given in Figure~3.

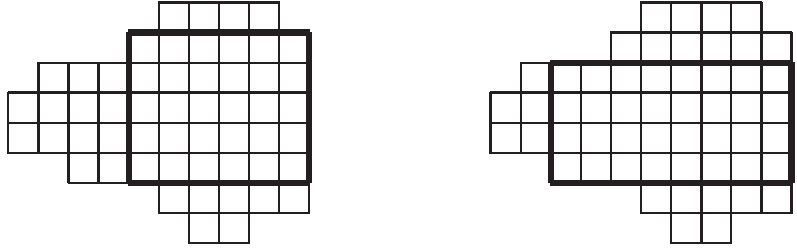
\begin{figure}[h!]
\begin{center}
 {\setlength{\unitlength}{0.8mm}
\begin{picture}(50,50)(0,0)
\put(30,0){\line(1,0){10}} \put(25,5){\line(1,0){25}}
\put(10,10){\line(1,0){40}} \put(0,15){\line(1,0){50}}
\put(0,20){\line(1,0){50}} \put(0,25){\line(1,0){50}}
\put(5,30){\line(1,0){45}} \put(20,35){\line(1,0){30}}
\put(15,30){\line(1,0){10}} \put(25,40){\line(1,0){20}}
\put(0,15){\line(0,1){10}}\put(5,15){\line(0,1){15}}
\put(10,10){\line(0,1){20}}\put(15,10){\line(0,1){20}}\put(20,10){\line(0,1){25}}
\put(25,5){\line(0,1){35}}\put(30,0){\line(0,1){40}}\put(35,0){\line(0,1){40}}\put(40,0){\line(0,1){40}}
\put(45,5){\line(0,1){35}} \put(50,5){\line(0,1){30}}
\linethickness{0.7mm}
\put(20,35){\line(1,0){30}}\put(20,10){\line(1,0){30}}
\put(20,10){\line(0,1){25}}\put(50,10){\line(0,1){25}}
\end{picture}}
\hspace{2cm} {\setlength{\unitlength}{0.8mm}
\begin{picture}(50,50)(0,0)
\put(30,0){\line(1,0){10}} \put(25,5){\line(1,0){25}}
\put(10,10){\line(1,0){40}} \put(0,15){\line(1,0){50}}
\put(0,20){\line(1,0){50}} \put(0,25){\line(1,0){50}}
\put(5,30){\line(1,0){45}} \put(20,35){\line(1,0){30}}
\put(15,30){\line(1,0){10}} \put(25,40){\line(1,0){20}}
\put(0,15){\line(0,1){10}}\put(5,15){\line(0,1){15}}
\put(10,10){\line(0,1){20}}\put(15,10){\line(0,1){20}}\put(20,10){\line(0,1){25}}
\put(25,5){\line(0,1){35}}\put(30,0){\line(0,1){40}}\put(35,0){\line(0,1){40}}\put(40,0){\line(0,1){40}}
\put(45,5){\line(0,1){35}} \put(50,5){\line(0,1){30}}
\linethickness{0.7mm}
\put(10,10){\line(1,0){40}}\put(10,30){\line(1,0){40}}
\put(10,10){\line(0,1){20}}\put(50,10){\line(0,1){20}}
\end{picture}}
\end{center}
\caption{\emph{left}: the $2$-th rectangle, \emph{right}: the $6$-th
rectangle.}
\end{figure}

\subsection{Coloring of fillings}

Let $F\in\N(T,{\bf m})$. The coloring of $F$ is the colored filling
obtained from $F$ by:
\begin{itemize}
\item coloring the cells of the empty columns,
\item for each $R_i\in Up(T)$ (resp., $R_i\in Low(T)$), coloring
the cells which are contained in the $i$-th rectangle and below
(resp., above) each $\textbf{1}$ in $R_i$. An illustration is given
in Figure~4.
\end{itemize}

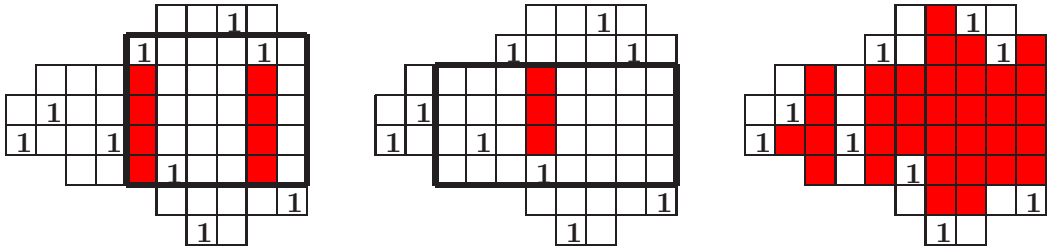
\begin{figure}[h!]\label{fig:coloring of a moon polyomino}
\begin{center}
 {\setlength{\unitlength}{0.8mm}
\begin{picture}(50,50)(0,0)
\put(20,10){\color{red}{\rule{4mm}{16mm}}}\put(40,10){\color{red}{\rule{4mm}{16mm}}}
\put(30,0){\line(1,0){10}} \put(25,5){\line(1,0){25}}
\put(10,10){\line(1,0){40}} \put(0,15){\line(1,0){50}}
\put(0,20){\line(1,0){50}} \put(0,25){\line(1,0){50}}
\put(5,30){\line(1,0){45}} \put(20,35){\line(1,0){30}}
\put(15,30){\line(1,0){10}} \put(25,40){\line(1,0){20}}
\put(0,15){\line(0,1){10}}\put(5,15){\line(0,1){15}}
\put(10,10){\line(0,1){20}}\put(15,10){\line(0,1){20}}\put(20,10){\line(0,1){25}}
\put(25,5){\line(0,1){35}}\put(30,0){\line(0,1){40}}\put(35,0){\line(0,1){40}}\put(40,0){\line(0,1){40}}
\put(45,5){\line(0,1){35}} \put(50,5){\line(0,1){30}}
\put(30,0){\makebox(6,4)[c]{\small \1}}
\put(45,5){\makebox(6,4)[c]{\small \1}}
\put(25,10){\makebox(6,4)[c]{\small \1}}
\put(15,15){\makebox(6,4)[c]{\small\1}}
\put(0,15){\makebox(6,4)[c]{\small \1}}
\put(5,20){\makebox(6,4)[c]{\small\1}}
\put(20,30){\makebox(6,4)[c]{\small\1}}
\put(40,30){\makebox(6,4)[c]{\small \1}}
\put(35,35){\makebox(6,4)[c]{\small \1}}
\linethickness{0.7mm}
\put(20,35){\line(1,0){30}}\put(20,10){\line(1,0){30}}
\put(20,10){\line(0,1){25}}\put(50,10){\line(0,1){25}}
\end{picture}}
\hspace{0.5cm} {\setlength{\unitlength}{0.8mm}
\begin{picture}(50,50)(0,0)
\put(25,15){\color{red}{\rule{4mm}{12mm}}}
\put(30,0){\line(1,0){10}} \put(25,5){\line(1,0){25}}
\put(10,10){\line(1,0){40}} \put(0,15){\line(1,0){50}}
\put(0,20){\line(1,0){50}} \put(0,25){\line(1,0){50}}
\put(5,30){\line(1,0){45}} \put(20,35){\line(1,0){30}}
\put(15,30){\line(1,0){10}} \put(25,40){\line(1,0){20}}
\put(0,15){\line(0,1){10}}\put(5,15){\line(0,1){15}}
\put(10,10){\line(0,1){20}}\put(15,10){\line(0,1){20}}\put(20,10){\line(0,1){25}}
\put(25,5){\line(0,1){35}}\put(30,0){\line(0,1){40}}\put(35,0){\line(0,1){40}}\put(40,0){\line(0,1){40}}
\put(45,5){\line(0,1){35}} \put(50,5){\line(0,1){30}}
\put(30,0){\makebox(6,4)[c]{\small \1}}
\put(45,5){\makebox(6,4)[c]{\small \1}}
\put(25,10){\makebox(6,4)[c]{\small \1}}
\put(15,15){\makebox(6,4)[c]{\small\1}}
\put(0,15){\makebox(6,4)[c]{\small \1}}
\put(5,20){\makebox(6,4)[c]{\small\1}}
\put(20,30){\makebox(6,4)[c]{\small\1}}
\put(40,30){\makebox(6,4)[c]{\small \1}}
\put(35,35){\makebox(6,4)[c]{\small \1}}
\linethickness{0.7mm}
\put(10,10){\line(1,0){40}}\put(10,30){\line(1,0){40}}
\put(10,10){\line(0,1){20}}\put(50,10){\line(0,1){20}}
\end{picture}}
\hspace{0.5cm}
 {\setlength{\unitlength}{0.8mm}
\begin{picture}(50,50)(0,0)
\put(30,5){\color{red}{\rule{4mm}{28mm}}}\put(35,5){\color{red}{\rule{4mm}{24mm}}}
\put(45,10){\color{red}{\rule{4mm}{20mm}}}\put(25,15){\color{red}{\rule{4mm}{12mm}}}
\put(5,15){\color{red}{\rule{4mm}{4mm}}}\put(10,10){\color{red}{\rule{4mm}{16mm}}}
\put(20,10){\color{red}{\rule{4mm}{16mm}}}\put(40,10){\color{red}{\rule{4mm}{16mm}}}
\put(30,0){\line(1,0){10}} \put(25,5){\line(1,0){25}}
\put(10,10){\line(1,0){40}} \put(0,15){\line(1,0){50}}
\put(0,20){\line(1,0){50}} \put(0,25){\line(1,0){50}}
\put(5,30){\line(1,0){45}} \put(20,35){\line(1,0){30}}
\put(15,30){\line(1,0){10}} \put(25,40){\line(1,0){20}}
\put(0,15){\line(0,1){10}}\put(5,15){\line(0,1){15}}
\put(10,10){\line(0,1){20}}\put(15,10){\line(0,1){20}}\put(20,10){\line(0,1){25}}
\put(25,5){\line(0,1){35}}\put(30,0){\line(0,1){40}}\put(35,0){\line(0,1){40}}\put(40,0){\line(0,1){40}}
\put(45,5){\line(0,1){35}} \put(50,5){\line(0,1){30}}
\put(30,0){\makebox(6,4)[c]{\small \1}}
\put(45,5){\makebox(6,4)[c]{\small \1}}
\put(25,10){\makebox(6,4)[c]{\small \1}}
\put(15,15){\makebox(6,4)[c]{\small\1}}
\put(0,15){\makebox(6,4)[c]{\small \1}}
\put(5,20){\makebox(6,4)[c]{\small\1}}
\put(20,30){\makebox(6,4)[c]{\small\1}}
\put(40,30){\makebox(6,4)[c]{\small \1}}
\put(35,35){\makebox(6,4)[c]{\small \1}}
\end{picture}}
\end{center}
\caption{\emph{left}: coloring induced by $R_2$, \emph{center}:
coloring induced by $R_6$, \emph{right}: full coloring.}
\end{figure}

The interest of coloring a $01$-filling is in the following result.
Let $c$ be a cell of $F$. If $c$ is filled with $1$ we denote by
$\luc(c;F)$ (resp., $\ruc(c;F)$) the numbers of uncolored cells to
the left (resp. right) and in the same row than the cell $c$ in the
coloring of $F$. If $c$ is empty, we set $\luc(c;F)=\ruc(c;F)=0$.
The proof of the following result is left to the reader.

\begin{prop}\label{prop:colored filling}
Let $F\in\N(T,{\bf m})$ and $c$ be a cell of $R_i$ filled with \1.
Then $\luc(c;F)$ (resp., $\ruc(c;F)$) is equal to
\begin{itemize}
\item if $R_i\in Up(T)$: the number of ascents (resp., descents) contained in the $i$-th rectangle of $F$ whose
North-east (resp., North-west) $1$ is in $c$,
\item if $R_i\in Low(T)$: the number of descents (resp. ascents)  contained in the $i$-th rectangle of $F$ whose South-east
(resp., South-west) $1$ is in $c$.
\end{itemize}
\end{prop}

It then follows from the above Proposition that
\begin{align}
\nec_2(F)=\sum_{c \,\in \,Up(F)} \luc(c;F)+\sum_{c\, \in \,Low(F)}\ruc(c;F),\label{eq:cro-sum}\\
\sec_2(F)=\sum_{c \,\in \,Up(F)}\ruc(c;F)+\sum_{c
\,\in\,Low(F)}\luc(c;F).\label{eq:ne-sum}
\end{align}

We can now describe our involution.
\subsection{The involution $\Phi$}

Let $F\in \N(T,{\bf m})$. We construct $\Phi(F)\in \N(T,{\bf m})$ by
the following process.  We start with the polyomino (empty filling)
$T$.

 (1) Color the columns of $T$ indexed by $\EC(F)$. We denote by $F'_0$ the
result.

(2) Suppose $R_{i_1}\prec R_{i_2} \prec \cdots \prec R_{i_s}$. For
$j$ from 1 to $s$, the (colored) filling $F'_{j}$ is obtained from
$F'_{j-1}$ by the following process:

\begin{itemize}
\item if $m_{i_j}=0$, then do nothing,
\item else, read the $m_{i_j}$ 1's in the ${i_j}$-th row of
$F$ from left to right and denote the number of uncolored cells
strictly
\begin{itemize}
\item to the left of the first~$1$ by $h_{0}$,
\item between the $u$-th 1 and the $(u+1)$-th 1, $1\leq u\leq m_{i_j}-1$, by $h_{u}$,
\item to the right of the last~$1$ by $h_{m_{i_j}}$.
\end{itemize}

 Then fill the $i_j\,$-th row of $F'_{j-1}$ with $m_{i_j}$ 1's in such a
way that the number of uncolored cells strictly
\begin{itemize}
\item to the left of the first~$1$ is $h_{m_{i_j}}$,
\item between the $u$-th 1 and the $(u+1)$-th 1, $1\leq u\leq m_{i_j}-1$, is $h_{m_{i_j}-u}$,
\item to the right of the last~$1$ is $h_{0}$.
\end{itemize}
Next, color the cells which are below (resp., above) the new $1$'s
and contained in the $i_j$-th rectangle if $R_{i_j}\in Up(T)$
(resp., $R_{i_j}\in Low(T)$).
\end{itemize}

(3) Set $\Phi(F)=F'_s$. For a better understanding, we give an
illustration. Suppose $F$ is the filling given below. \begin{center}
{\setlength{\unitlength}{0.8mm}
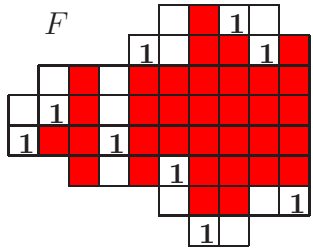
\begin{picture}(50,40)(0,0)
\put(5,35){\makebox(6,4)[c]{$F$}}
\put(30,5){\color{red}{\rule{4mm}{28mm}}}\put(35,5){\color{red}{\rule{4mm}{24mm}}}
\put(45,10){\color{red}{\rule{4mm}{20mm}}}\put(25,15){\color{red}{\rule{4mm}{12mm}}}
\put(5,15){\color{red}{\rule{4mm}{4mm}}}\put(10,10){\color{red}{\rule{4mm}{16mm}}}
\put(20,10){\color{red}{\rule{4mm}{16mm}}}\put(40,10){\color{red}{\rule{4mm}{16mm}}}
\put(30,0){\line(1,0){10}} \put(25,5){\line(1,0){25}}
\put(10,10){\line(1,0){40}} \put(0,15){\line(1,0){50}}
\put(0,20){\line(1,0){50}} \put(0,25){\line(1,0){50}}
\put(5,30){\line(1,0){45}} \put(20,35){\line(1,0){30}}
\put(15,30){\line(1,0){10}} \put(25,40){\line(1,0){20}}
\put(0,15){\line(0,1){10}}\put(5,15){\line(0,1){15}}
\put(10,10){\line(0,1){20}}\put(15,10){\line(0,1){20}}\put(20,10){\line(0,1){25}}
\put(25,5){\line(0,1){35}}\put(30,0){\line(0,1){40}}\put(35,0){\line(0,1){40}}\put(40,0){\line(0,1){40}}
\put(45,5){\line(0,1){35}} \put(50,5){\line(0,1){30}}
\put(30,0){\makebox(6,4)[c]{\small \1}}
\put(45,5){\makebox(6,4)[c]{\small \1}}
\put(25,10){\makebox(6,4)[c]{\small \1}}
\put(15,15){\makebox(6,4)[c]{\small\1}}
\put(0,15){\makebox(6,4)[c]{\small \1}}
\put(5,20){\makebox(6,4)[c]{\small\1}}
\put(20,30){\makebox(6,4)[c]{\small\1}}
\put(40,30){\makebox(6,4)[c]{\small \1}}
\put(35,35){\makebox(6,4)[c]{\small \1}}
\end{picture}}\end{center}

Then the step-by-step construction of $\Phi(F)$ goes as follows.

\begin{center}
 {\setlength{\unitlength}{0.8mm}
\begin{picture}(50,45)(0,0)
\put(5,35){\makebox(6,4)[c]{$F'_0$}}
\put(10,10){\color{red}{\rule{4mm}{16mm}}}
\put(30,0){\line(1,0){10}} \put(25,5){\line(1,0){25}}
\put(10,10){\line(1,0){40}} \put(0,15){\line(1,0){50}}
\put(0,20){\line(1,0){50}} \put(0,25){\line(1,0){50}}
\put(5,30){\line(1,0){45}} \put(20,35){\line(1,0){30}}
\put(15,30){\line(1,0){10}} \put(25,40){\line(1,0){20}}
\put(0,15){\line(0,1){10}}\put(5,15){\line(0,1){15}}
\put(10,10){\line(0,1){20}}\put(15,10){\line(0,1){20}}\put(20,10){\line(0,1){25}}
\put(25,5){\line(0,1){35}}\put(30,0){\line(0,1){40}}\put(35,0){\line(0,1){40}}\put(40,0){\line(0,1){40}}
\put(45,5){\line(0,1){35}} \put(50,5){\line(0,1){30}}
\end{picture}}
\hspace{1cm}
 {\setlength{\unitlength}{0.8mm}
\begin{picture}(50,45)(0,0)
\put(5,35){\makebox(6,4)[c]{$F'_1$}}
\put(10,10){\color{red}{\rule{4mm}{16mm}}}
\put(35,5){\color{red}{\rule{4mm}{28mm}}}
\put(30,0){\line(1,0){10}} \put(25,5){\line(1,0){25}}
\put(10,10){\line(1,0){40}} \put(0,15){\line(1,0){50}}
\put(0,20){\line(1,0){50}} \put(0,25){\line(1,0){50}}
\put(5,30){\line(1,0){45}} \put(20,35){\line(1,0){30}}
\put(15,30){\line(1,0){10}} \put(25,40){\line(1,0){20}}
\put(0,15){\line(0,1){10}}\put(5,15){\line(0,1){15}}
\put(10,10){\line(0,1){20}}\put(15,10){\line(0,1){20}}\put(20,10){\line(0,1){25}}
\put(25,5){\line(0,1){35}}\put(30,0){\line(0,1){40}}\put(35,0){\line(0,1){40}}\put(40,0){\line(0,1){40}}
\put(45,5){\line(0,1){35}} \put(50,5){\line(0,1){30}}
\put(35,0){\makebox(6,4)[c]{\small \1}}
\linethickness{0.7mm}
\put(30,0){\line(1,0){10}}\put(30,40){\line(1,0){10}}
\put(30,0){\line(0,1){40}}\put(40,0){\line(0,1){40}}
\end{picture}}
\hspace{1cm}
 {\setlength{\unitlength}{0.8mm}
\begin{picture}(50,45)(0,0)
\put(5,35){\makebox(6,4)[c]{$F'_2$}}
\put(10,10){\color{red}{\rule{4mm}{16mm}}}
\put(35,5){\color{red}{\rule{4mm}{28mm}}}
\put(30,5){\color{red}{\rule{4mm}{24mm}}}
\put(30,0){\line(1,0){10}} \put(25,5){\line(1,0){25}}
\put(10,10){\line(1,0){40}} \put(0,15){\line(1,0){50}}
\put(0,20){\line(1,0){50}} \put(0,25){\line(1,0){50}}
\put(5,30){\line(1,0){45}} \put(20,35){\line(1,0){30}}
\put(15,30){\line(1,0){10}} \put(25,40){\line(1,0){20}}
\put(0,15){\line(0,1){10}}\put(5,15){\line(0,1){15}}
\put(10,10){\line(0,1){20}}\put(15,10){\line(0,1){20}}\put(20,10){\line(0,1){25}}
\put(25,5){\line(0,1){35}}\put(30,0){\line(0,1){40}}\put(35,0){\line(0,1){40}}\put(40,0){\line(0,1){40}}
\put(45,5){\line(0,1){35}} \put(50,5){\line(0,1){30}}
\put(35,0){\makebox(6,4)[c]{\small \1}}
\put(30,35){\makebox(6,4)[c]{\small \1}}
\linethickness{0.7mm}
\put(25,5){\line(0,1){35}}\put(45,5){\line(0,1){35}}
\put(25,5){\line(1,0){20}}\put(25,40){\line(1,0){20}}
\end{picture}}
\end{center}

\begin{center}
 {\setlength{\unitlength}{0.8mm}
\begin{picture}(50,50)(0,0)
\put(5,35){\makebox(6,4)[c]{$F'_3$}}
\put(10,10){\color{red}{\rule{4mm}{16mm}}}
\put(35,5){\color{red}{\rule{4mm}{28mm}}}
\put(30,5){\color{red}{\rule{4mm}{24mm}}}
\put(25,10){\color{red}{\rule{4mm}{20mm}}}
\put(30,0){\line(1,0){10}} \put(25,5){\line(1,0){25}}
\put(10,10){\line(1,0){40}} \put(0,15){\line(1,0){50}}
\put(0,20){\line(1,0){50}} \put(0,25){\line(1,0){50}}
\put(5,30){\line(1,0){45}} \put(20,35){\line(1,0){30}}
\put(15,30){\line(1,0){10}} \put(25,40){\line(1,0){20}}
\put(0,15){\line(0,1){10}}\put(5,15){\line(0,1){15}}
\put(10,10){\line(0,1){20}}\put(15,10){\line(0,1){20}}\put(20,10){\line(0,1){25}}
\put(25,5){\line(0,1){35}}\put(30,0){\line(0,1){40}}\put(35,0){\line(0,1){40}}\put(40,0){\line(0,1){40}}
\put(45,5){\line(0,1){35}} \put(50,5){\line(0,1){30}}
\put(35,0){\makebox(6,4)[c]{\small \1}}
\put(30,35){\makebox(6,4)[c]{\small \1}}
\put(25,5){\makebox(6,4)[c]{\small \1}}
\linethickness{0.7mm}
\put(25,5){\line(0,1){30}}\put(50,5){\line(0,1){30}}
\put(25,5){\line(1,0){25}}\put(25,35){\line(1,0){25}}
\end{picture}}
\hspace{1cm}
 {\setlength{\unitlength}{0.8mm}
\begin{picture}(50,50)(0,0)
\put(5,35){\makebox(6,4)[c]{$F'_4$}}
\put(10,10){\color{red}{\rule{4mm}{16mm}}}
\put(35,5){\color{red}{\rule{4mm}{28mm}}}
\put(30,5){\color{red}{\rule{4mm}{24mm}}}
\put(25,10){\color{red}{\rule{4mm}{20mm}}}
\put(20,10){\color{red}{\rule{4mm}{16mm}}}
\put(45,10){\color{red}{\rule{4mm}{16mm}}}
\put(30,0){\line(1,0){10}} \put(25,5){\line(1,0){25}}
\put(10,10){\line(1,0){40}} \put(0,15){\line(1,0){50}}
\put(0,20){\line(1,0){50}} \put(0,25){\line(1,0){50}}
\put(5,30){\line(1,0){45}} \put(20,35){\line(1,0){30}}
\put(15,30){\line(1,0){10}} \put(25,40){\line(1,0){20}}
\put(0,15){\line(0,1){10}}\put(5,15){\line(0,1){15}}
\put(10,10){\line(0,1){20}}\put(15,10){\line(0,1){20}}\put(20,10){\line(0,1){25}}
\put(25,5){\line(0,1){35}}\put(30,0){\line(0,1){40}}\put(35,0){\line(0,1){40}}\put(40,0){\line(0,1){40}}
\put(45,5){\line(0,1){35}} \put(50,5){\line(0,1){30}}
\put(35,0){\makebox(6,4)[c]{\small \1}}
\put(30,35){\makebox(6,4)[c]{\small \1}}
\put(25,5){\makebox(6,4)[c]{\small \1}}
\put(20,30){\makebox(6,4)[c]{\small \1}}
\put(45,30){\makebox(6,4)[c]{\small \1}}
\linethickness{0.7mm}
\put(20,10){\line(0,1){25}}\put(50,10){\line(0,1){25}}
\put(20,10){\line(1,0){30}}\put(20,35){\line(1,0){30}}
\end{picture}}
\hspace{1cm}
{\setlength{\unitlength}{0.8mm}
\begin{picture}(50,50)(0,0)
\put(5,35){\makebox(6,4)[c]{$F'_5$}}
\put(10,10){\color{red}{\rule{4mm}{16mm}}}
\put(35,5){\color{red}{\rule{4mm}{28mm}}}
\put(30,5){\color{red}{\rule{4mm}{24mm}}}
\put(25,10){\color{red}{\rule{4mm}{20mm}}}
\put(20,10){\color{red}{\rule{4mm}{16mm}}}
\put(45,10){\color{red}{\rule{4mm}{16mm}}}
\put(15,15){\color{red}{\rule{4mm}{12mm}}}
\put(30,0){\line(1,0){10}} \put(25,5){\line(1,0){25}}
\put(10,10){\line(1,0){40}} \put(0,15){\line(1,0){50}}
\put(0,20){\line(1,0){50}} \put(0,25){\line(1,0){50}}
\put(5,30){\line(1,0){45}} \put(20,35){\line(1,0){30}}
\put(15,30){\line(1,0){10}} \put(25,40){\line(1,0){20}}
\put(0,15){\line(0,1){10}}\put(5,15){\line(0,1){15}}
\put(10,10){\line(0,1){20}}\put(15,10){\line(0,1){20}}\put(20,10){\line(0,1){25}}
\put(25,5){\line(0,1){35}}\put(30,0){\line(0,1){40}}\put(35,0){\line(0,1){40}}\put(40,0){\line(0,1){40}}
\put(45,5){\line(0,1){35}} \put(50,5){\line(0,1){30}}
\put(35,0){\makebox(6,4)[c]{\small \1}}
\put(30,35){\makebox(6,4)[c]{\small \1}}
\put(25,5){\makebox(6,4)[c]{\small \1}}
\put(20,30){\makebox(6,4)[c]{\small \1}}
\put(45,30){\makebox(6,4)[c]{\small \1}}
\put(15,10){\makebox(6,4)[c]{\small \1}}
\linethickness{0.7mm}
\put(10,10){\line(0,1){20}}\put(50,10){\line(0,1){20}}
\put(10,10){\line(1,0){40}}\put(10,30){\line(1,0){40}}
\end{picture}}
\end{center}

\begin{figure}[h!]\label{fig:invo-fillings}

\vspace{-0.1cm}

\vspace{-0.3cm}

\vspace{-0.3cm}
\begin{center}
 {\setlength{\unitlength}{0.8mm}
\begin{picture}(50,50)(0,0)
\put(5,35){\makebox(6,4)[c]{$F'_6$}}
\put(10,10){\color{red}{\rule{4mm}{16mm}}}
\put(35,5){\color{red}{\rule{4mm}{28mm}}}
\put(30,5){\color{red}{\rule{4mm}{24mm}}}
\put(25,10){\color{red}{\rule{4mm}{20mm}}}
\put(20,10){\color{red}{\rule{4mm}{16mm}}}
\put(45,10){\color{red}{\rule{4mm}{16mm}}}
\put(15,15){\color{red}{\rule{4mm}{12mm}}}
\put(30,0){\line(1,0){10}} \put(25,5){\line(1,0){25}}
\put(10,10){\line(1,0){40}} \put(0,15){\line(1,0){50}}
\put(0,20){\line(1,0){50}} \put(0,25){\line(1,0){50}}
\put(5,30){\line(1,0){45}} \put(20,35){\line(1,0){30}}
\put(15,30){\line(1,0){10}} \put(25,40){\line(1,0){20}}
\put(0,15){\line(0,1){10}}\put(5,15){\line(0,1){15}}
\put(10,10){\line(0,1){20}}\put(15,10){\line(0,1){20}}\put(20,10){\line(0,1){25}}
\put(25,5){\line(0,1){35}}\put(30,0){\line(0,1){40}}\put(35,0){\line(0,1){40}}\put(40,0){\line(0,1){40}}
\put(45,5){\line(0,1){35}} \put(50,5){\line(0,1){30}}
\put(35,0){\makebox(6,4)[c]{\small \1}}
\put(30,35){\makebox(6,4)[c]{\small \1}}
\put(25,5){\makebox(6,4)[c]{\small \1}}
\put(20,30){\makebox(6,4)[c]{\small \1}}
\put(45,30){\makebox(6,4)[c]{\small \1}}
\put(15,10){\makebox(6,4)[c]{\small \1}}
\linethickness{0.7mm}
\put(5,15){\line(0,1){15}}\put(50,15){\line(0,1){15}}
\put(5,15){\line(1,0){45}}\put(5,30){\line(1,0){45}}
\end{picture}}
\hspace{1cm}
{\setlength{\unitlength}{0.8mm}
\begin{picture}(50,50)(0,0)
\put(5,35){\makebox(6,4)[c]{$F'_7$}}
\put(10,10){\color{red}{\rule{4mm}{16mm}}}
\put(35,5){\color{red}{\rule{4mm}{28mm}}}
\put(30,5){\color{red}{\rule{4mm}{24mm}}}
\put(25,10){\color{red}{\rule{4mm}{20mm}}}
\put(20,10){\color{red}{\rule{4mm}{16mm}}}
\put(45,10){\color{red}{\rule{4mm}{16mm}}}
\put(15,15){\color{red}{\rule{4mm}{12mm}}}
\put(5,15){\color{red}{\rule{4mm}{4mm}}}
\put(30,0){\line(1,0){10}} \put(25,5){\line(1,0){25}}
\put(10,10){\line(1,0){40}} \put(0,15){\line(1,0){50}}
\put(0,20){\line(1,0){50}} \put(0,25){\line(1,0){50}}
\put(5,30){\line(1,0){45}} \put(20,35){\line(1,0){30}}
\put(15,30){\line(1,0){10}} \put(25,40){\line(1,0){20}}
\put(0,15){\line(0,1){10}}\put(5,15){\line(0,1){15}}
\put(10,10){\line(0,1){20}}\put(15,10){\line(0,1){20}}\put(20,10){\line(0,1){25}}
\put(25,5){\line(0,1){35}}\put(30,0){\line(0,1){40}}\put(35,0){\line(0,1){40}}\put(40,0){\line(0,1){40}}
\put(45,5){\line(0,1){35}} \put(50,5){\line(0,1){30}}
\put(35,0){\makebox(6,4)[c]{\small \1}}
\put(30,35){\makebox(6,4)[c]{\small \1}}
\put(25,5){\makebox(6,4)[c]{\small \1}}
\put(20,30){\makebox(6,4)[c]{\small \1}}
\put(45,30){\makebox(6,4)[c]{\small \1}}
\put(15,10){\makebox(6,4)[c]{\small \1}}
\put(5,20){\makebox(6,4)[c]{\small \1}}
\linethickness{0.7mm}
\put(0,15){\line(0,1){10}}\put(50,15){\line(0,1){10}}
\put(0,15){\line(1,0){50}}\put(0,25){\line(1,0){50}}
\end{picture}}
\hspace{1cm}
{\setlength{\unitlength}{0.8mm}
\begin{picture}(50,50)(0,0)
\put(2,35){\makebox(6,4)[c]{$F'_8=\Phi(F)$}}
\put(10,10){\color{red}{\rule{4mm}{16mm}}}
\put(35,5){\color{red}{\rule{4mm}{28mm}}}
\put(30,5){\color{red}{\rule{4mm}{24mm}}}
\put(25,10){\color{red}{\rule{4mm}{20mm}}}
\put(20,10){\color{red}{\rule{4mm}{16mm}}}
\put(45,10){\color{red}{\rule{4mm}{16mm}}}
\put(15,15){\color{red}{\rule{4mm}{12mm}}}
\put(5,15){\color{red}{\rule{4mm}{4mm}}}
\put(30,0){\line(1,0){10}} \put(25,5){\line(1,0){25}}
\put(10,10){\line(1,0){40}} \put(0,15){\line(1,0){50}}
\put(0,20){\line(1,0){50}} \put(0,25){\line(1,0){50}}
\put(5,30){\line(1,0){45}} \put(20,35){\line(1,0){30}}
\put(15,30){\line(1,0){10}} \put(25,40){\line(1,0){20}}
\put(0,15){\line(0,1){10}}\put(5,15){\line(0,1){15}}
\put(10,10){\line(0,1){20}}\put(15,10){\line(0,1){20}}\put(20,10){\line(0,1){25}}
\put(25,5){\line(0,1){35}}\put(30,0){\line(0,1){40}}\put(35,0){\line(0,1){40}}\put(40,0){\line(0,1){40}}
\put(45,5){\line(0,1){35}} \put(50,5){\line(0,1){30}}
\put(35,0){\makebox(6,4)[c]{\small \1}}
\put(30,35){\makebox(6,4)[c]{\small \1}}
\put(25,5){\makebox(6,4)[c]{\small \1}}
\put(20,30){\makebox(6,4)[c]{\small \1}}
\put(45,30){\makebox(6,4)[c]{\small \1}}
\put(15,10){\makebox(6,4)[c]{\small \1}}
\put(5,20){\makebox(6,4)[c]{\small \1}}
\put(0,15){\makebox(6,4)[c]{\small \1}}
\put(40,15){\makebox(6,4)[c]{\small \1}}
\linethickness{0.7mm}
\put(0,15){\line(0,1){5}}\put(50,15){\line(0,1){5}}
\put(0,15){\line(1,0){50}}\put(0,20){\line(1,0){50}}
\end{picture}}
\end{center}
\caption{ The step-by-step construction of $\Phi(F)$}
\end{figure}
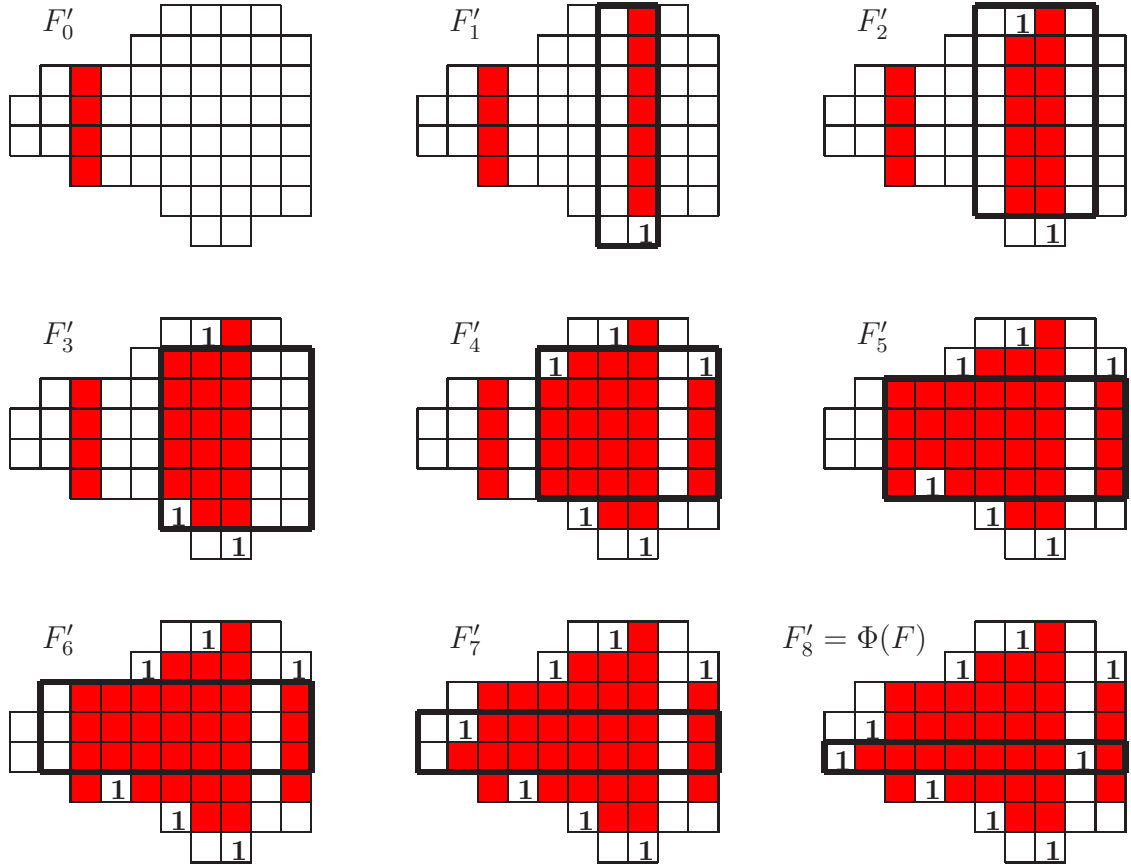

It is easy to check that the map $\Phi$ is an involution on
$\N(T,{\bf m})$ and $\EC(\Phi(F))=\EC(F)$ for any $F$. Invoking
Equations \eqref{eq:cro-sum} and \eqref{eq:ne-sum} lead to the
following result which imply Theorem~\ref{thm:refinement main1}.

\begin{prop}\label{prop:Phi}
For any moon polyomino $T$, the map $\Phi$ is an involution on
$\N(T,{\bf m})$ such that for any $F\in\N(T,{\bf m})$  we have
\begin{align*}
\EC(\Phi(F))=\EC(F)\;,\;\nec_2(\Phi(F))=\sec_2(F)\;,\;\sec_2(\Phi(F))=\nec_2(F).
\end{align*}
\end{prop}

One can also give an alternating proof of
Theorem~\ref{thm:refinement main1} by computing the joint
distribution of $(\nec_2,\sec_2)$. This is the content of the
following section.

\section{Distribution of ascents and descents over $\N(T,{\bf m},A)$}

For nonnegative integers $n$ and $k$, let ${n\brack k}_{p,q}$ be the
$p,q$-Gaussian coefficient defined by
$$
{n\brack k}_{p,q}=\frac{[n]_{p,q}!}{[k]_{p,q}!\;[n-k]_{p,q}!},
$$
where, as usual in $p,q$-theory, the $p,q$-integer $[r]_{p,q}$ is
given by
$$[r]_{p,q}:=\frac{p^i-q^i}{p-q}=(p^{i-1}+p^{i-2}q+\cdots+p^jq^{i-j-1}+\cdots+pq^{i-2}+q^{i-1}),$$
and the $p,q$-factorial $[r]_{p,q}!$  by
$[r]_{p,q}!:=\prod_{i=1}^r[i]_{p,q}$.\\

 Let $T$ be a moon polyomino
with $s$ rows and $t$ columns, ${\bf m}=(m_1,\ldots,m_s)$ a $s$-uple
of nonnegative integers and $A$ a subset of $[t]$. Suppose
$R_{i_1}\prec R_{i_2} \prec \cdots \prec R_{i_s}$. Then for
$j\in[s]$, define $h_{i_j}$ by
\begin{equation}\label{eq:h_i}
h_{i_j}=r_{i_j}-(m_{i_1}+m_{i_2}+\cdots+m_{i_{j-1}})-a_{{i_j}},
\end{equation}
where  $r_{i_j}$ is the length of row $R_{i_j}$ and $a_{{i_j}}$ is
the number of indices $k\in A$ such that the column labeled $k$
intersect the row labeled $i_j$. Then we have the following result.

\begin{thm}\label{thm:distribution}
For any moon polyomino $T$, the distribution of the joint statistic
$(\nec_2,\sec_2)$ over $\N(T,{\bf m};A)$ is given by
\begin{equation}\label{eq:symfilling}
\sum_{F\in\,\N(T,{\bf
m};A)}p^{\nec_2(F)}q^{\sec_2(F)}=\prod_{j=1}^s{h_j\brack m_j}_{p,q},
\end{equation}
where $h_j$ is defined by \eqref{eq:h_i}. In particular,
\begin{equation}
|\N(T,{\bf m};A)|=\prod_{j=1}^s{h_j\choose m_j},
\end{equation}
\end{thm}

For instance, suppose $T$ be the moon polyominoe given below and
$A=\{2\}$.
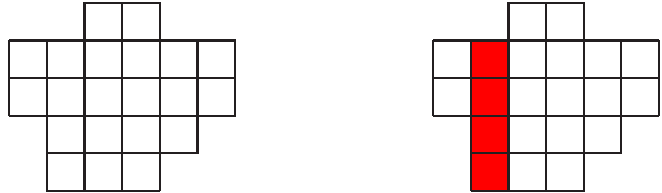
\begin{center}
 {\setlength{\unitlength}{1mm}
\begin{picture}(35,28)(0,0)
\put(5,0){\line(1,0){15}}\put(5,5){\line(1,0){20}}\put(0,10){\line(1,0){30}}
\put(0,15){\line(1,0){30}}\put(0,20){\line(1,0){30}}\put(10,25){\line(1,0){10}}
\put(0,10){\line(0,1){10}}\put(5,0){\line(0,1){20}}\put(10,0){\line(0,1){25}}
\put(15,0){\line(0,1){25}}\put(20,0){\line(0,1){25}}\put(25,5){\line(0,1){15}}
\put(30,10){\line(0,1){10}}
\end{picture}}\hspace{2cm}{\setlength{\unitlength}{1mm}
\begin{picture}(35,25)(0,0)
\put(5,0){\color{red}{\rule{5mm}{20mm}}}
\put(5,0){\line(1,0){15}}\put(5,5){\line(1,0){20}}\put(0,10){\line(1,0){30}}
\put(0,15){\line(1,0){30}}\put(0,20){\line(1,0){30}}\put(10,25){\line(1,0){10}}
\put(0,10){\line(0,1){10}}\put(5,0){\line(0,1){20}}\put(10,0){\line(0,1){25}}
\put(15,0){\line(0,1){25}}\put(20,0){\line(0,1){25}}\put(25,5){\line(0,1){15}}
\put(30,10){\line(0,1){10}}
\end{picture}}
\end{center}
Then we have:
\begin{itemize}
\item $R_1\prec R_5\prec R_4\prec R_2 \prec R_3$ thus  $i_1=1$, $i_2=5$, $i_3=4$,
$i_4=2$, $i_5=3$.
\item The column labeled $2$ intersect the rows labeled $2,3,4,5$, thus $a_1=0$,
$a_2=a_3=a_4=a_5=1$.
\end{itemize}
Suppose ${\bf m}=(1,2,1,0,1)$. We then have
\begin{align*}
h_{i_1}=h_1&=r_{1}-a_{1}=2,\\
h_{i_2}=h_5&=r_{5}-m_{1}-a_{5}=1,\\
h_{i_3}=h_4&=r_{4}-(m_{1}+m_{5})-a_{4}=1,\\
h_{i_4}=h_2&=r_{2}-(m_{1}+m_{5}+m_{4})-a_{2}=3,\\
h_{i_5}=h_3&=r_{3}-(m_{1}+m_{5}+m_{4}+m_{2})-a_{3}=1.
\end{align*}
It follows that
$$\prod_{j=1}^5{h_j\brack m_j}_{p,q}={2\brack
1}_{p,q}{3\brack 2}_{p,q}{1\brack 1}_{p,q}{1\brack 0}_{p,q}{1\brack
1}_{p,q}=p^3+2p^2q+2pq^2+q^3.$$

On the other hand, the fillings in $\N(T,{\bf m}, A)$ and the
corresponding values of $\nec_2$ and $\sec_2$ are listed below.

\begin{center}
 {\setlength{\unitlength}{1mm}
\begin{picture}(35,25)(0,0)
\put(5,0){\color{red}{\rule{5mm}{20mm}}}
\put(5,0){\line(1,0){15}}\put(5,5){\line(1,0){20}}\put(0,10){\line(1,0){30}}
\put(0,15){\line(1,0){30}}\put(0,20){\line(1,0){30}}\put(10,25){\line(1,0){10}}
\put(0,10){\line(0,1){10}}\put(5,0){\line(0,1){20}}\put(10,0){\line(0,1){25}}
\put(15,0){\line(0,1){25}}\put(20,0){\line(0,1){25}}\put(25,5){\line(0,1){15}}
\put(30,10){\line(0,1){10}}
\put(10,20){\makebox(6,4)[c]{\small\1}}
\put(0,15){\makebox(6,4)[c]{\small\1}}
\put(20,15){\makebox(6,4)[c]{\small \1}}
\put(25,10){\makebox(6,4)[c]{\small \1}}
\put(15,0){\makebox(6,4)[c]{\small \1}}
\put(5,-7){\makebox(15,4)[c]{ $\nec_2=0\;,\;\sec_2=3$}}
\end{picture}}
\hspace{1cm}
  {\setlength{\unitlength}{1mm}
\begin{picture}(35,25)(0,0)
\put(5,0){\color{red}{\rule{5mm}{20mm}}}
\put(5,0){\line(1,0){15}}\put(5,5){\line(1,0){20}}\put(0,10){\line(1,0){30}}
\put(0,15){\line(1,0){30}}\put(0,20){\line(1,0){30}}\put(10,25){\line(1,0){10}}
\put(0,10){\line(0,1){10}}\put(5,0){\line(0,1){20}}\put(10,0){\line(0,1){25}}
\put(15,0){\line(0,1){25}}\put(20,0){\line(0,1){25}}\put(25,5){\line(0,1){15}}
\put(30,10){\line(0,1){10}}
\put(10,20){\makebox(6,4)[c]{\small\1}}
\put(0,15){\makebox(6,4)[c]{\small\1}}
\put(25,15){\makebox(6,4)[c]{\small \1}}
\put(20,10){\makebox(6,4)[c]{\small \1}}
\put(15,0){\makebox(6,4)[c]{\small \1}}
\put(5,-7){\makebox(15,4)[c]{ $\nec_2=1\;,\;\sec_2=2$}}
\end{picture}}
\hspace{1cm}
  {\setlength{\unitlength}{1mm}
\begin{picture}(35,25)(0,0)
\put(5,0){\color{red}{\rule{5mm}{20mm}}}
\put(5,0){\line(1,0){15}}\put(5,5){\line(1,0){20}}\put(0,10){\line(1,0){30}}
\put(0,15){\line(1,0){30}}\put(0,20){\line(1,0){30}}\put(10,25){\line(1,0){10}}
\put(0,10){\line(0,1){10}}\put(5,0){\line(0,1){20}}\put(10,0){\line(0,1){25}}
\put(15,0){\line(0,1){25}}\put(20,0){\line(0,1){25}}\put(25,5){\line(0,1){15}}
\put(30,10){\line(0,1){10}}
\put(10,20){\makebox(6,4)[c]{\small\1}}
\put(20,15){\makebox(6,4)[c]{\small\1}}
\put(25,15){\makebox(6,4)[c]{\small \1}}
\put(0,10){\makebox(6,4)[c]{\small \1}}
\put(15,0){\makebox(6,4)[c]{\small \1}}
\put(5,-7){\makebox(15,4)[c]{ $\nec_2=2\;,\;\sec_2=1$}}
\end{picture}}
\end{center}
\vspace{1cm}
\begin{center}
 {\setlength{\unitlength}{1mm}
\begin{picture}(35,25)(0,0)
\put(5,0){\color{red}{\rule{5mm}{20mm}}}
\put(5,0){\line(1,0){15}}\put(5,5){\line(1,0){20}}\put(0,10){\line(1,0){30}}
\put(0,15){\line(1,0){30}}\put(0,20){\line(1,0){30}}\put(10,25){\line(1,0){10}}
\put(0,10){\line(0,1){10}}\put(5,0){\line(0,1){20}}\put(10,0){\line(0,1){25}}
\put(15,0){\line(0,1){25}}\put(20,0){\line(0,1){25}}\put(25,5){\line(0,1){15}}
\put(30,10){\line(0,1){10}}
\put(15,20){\makebox(6,4)[c]{\small\1}}
\put(0,15){\makebox(6,4)[c]{\small\1}}
\put(20,15){\makebox(6,4)[c]{\small \1}}
\put(25,10){\makebox(6,4)[c]{\small \1}}
\put(10,0){\makebox(6,4)[c]{\small \1}}
\put(5,-7){\makebox(15,4)[c]{ $\nec_2=1\;,\;\sec_2=2$}}
\end{picture}}
\hspace{1cm}
  {\setlength{\unitlength}{1mm}
\begin{picture}(35,25)(0,0)
\put(5,0){\color{red}{\rule{5mm}{20mm}}}
\put(5,0){\line(1,0){15}}\put(5,5){\line(1,0){20}}\put(0,10){\line(1,0){30}}
\put(0,15){\line(1,0){30}}\put(0,20){\line(1,0){30}}\put(10,25){\line(1,0){10}}
\put(0,10){\line(0,1){10}}\put(5,0){\line(0,1){20}}\put(10,0){\line(0,1){25}}
\put(15,0){\line(0,1){25}}\put(20,0){\line(0,1){25}}\put(25,5){\line(0,1){15}}
\put(30,10){\line(0,1){10}}
\put(15,20){\makebox(6,4)[c]{\small\1}}
\put(0,15){\makebox(6,4)[c]{\small\1}}
\put(25,15){\makebox(6,4)[c]{\small \1}}
\put(20,10){\makebox(6,4)[c]{\small \1}}
\put(10,0){\makebox(6,4)[c]{\small \1}}
\put(5,-7){\makebox(15,4)[c]{ $\nec_2=2\;,\;\sec_2=1$}}
\end{picture}}
\hspace{1cm}
  {\setlength{\unitlength}{1mm}
\begin{picture}(35,25)(0,0)
\put(5,0){\color{red}{\rule{5mm}{20mm}}}
\put(5,0){\line(1,0){15}}\put(5,5){\line(1,0){20}}\put(0,10){\line(1,0){30}}
\put(0,15){\line(1,0){30}}\put(0,20){\line(1,0){30}}\put(10,25){\line(1,0){10}}
\put(0,10){\line(0,1){10}}\put(5,0){\line(0,1){20}}\put(10,0){\line(0,1){25}}
\put(15,0){\line(0,1){25}}\put(20,0){\line(0,1){25}}\put(25,5){\line(0,1){15}}
\put(30,10){\line(0,1){10}}
\put(15,20){\makebox(6,4)[c]{\small\1}}
\put(20,15){\makebox(6,4)[c]{\small\1}}
\put(25,15){\makebox(6,4)[c]{\small \1}}
\put(0,10){\makebox(6,4)[c]{\small \1}}
\put(10,0){\makebox(6,4)[c]{\small \1}}
\put(5,-7){\makebox(15,4)[c]{ $\nec_2=3\;,\;\sec_2=0$}}
\end{picture}}
\end{center}
\vspace{1cm}

Summing up we get $\sum_{F\in\N(T,{\bf m},
A)}p^{\nec_{2}(F)}q^{\sec_{2}(F)}=p^3+2p^2q+2pq^2+q^3$, as
desired.\\

We don't give a rigorous prove of Theorem~\ref{thm:distribution} but
just a sketch.

\subsection*{Sketch of the proof of Theorem~\ref{thm:distribution}} If $n$ and $k$ are positive
integers, we will denote by $\C_k(n)$ the set of compositions of $n$
into $k$ nonnegative parts. Recall that a element in $\C_k(n)$ is
just a $k$-uple $(b_1,b_2,\ldots,b_k)$ of nonnegative integers such
that $b_1+b_2+\cdots+b_k=n$. We will construct a bijection
$$f:\N(T,{\bf m};A)\to
\C_{m_{1}+1}(h_{1})\times\C_{m_{2}+1}(h_{2})\times \cdots\times
\C_{m_{s}+1}(h_{s})$$ which keeps track of the statistics $\nec_2$
and $\sec_2$. To each $F\in\N(T,{\bf m};A)$ we associate the
sequence of compositions $(c^{(1)},c^{(2)},\ldots,c^{(s)})$, where
for $i=1,\ldots,s$, $c^i=(c^{(i)}_1,\ldots,c^{(i)}_{m_i+1})$ is
defined by $c^{(i)}_1$ (resp., $c^{(i)}_j$ for $j=2\ldots m_i$,
$c^{(i)}_{m_i+1}$) is the number of uncolored cells to the left of
the first~$1$ (resp., between the $j$-th 1 and the $(j+1)$-th 1, to
the right of the last~$1$) of $R_i$, the ${i}$-th row of the
coloring of~$F$. An example is given below.

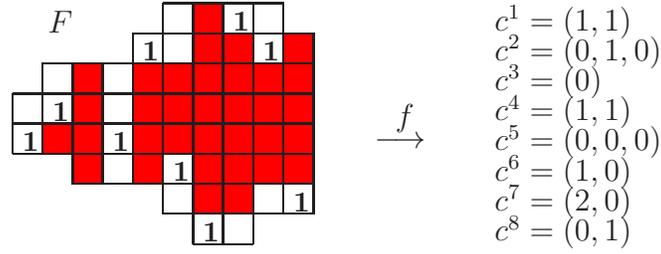
\begin{figure}[h!]
\begin{center}
 {\setlength{\unitlength}{0.8mm}
\begin{picture}(50,50)(0,-5)
\put(5,35){\makebox(6,4)[c]{$F$}}
\put(30,5){\color{red}{\rule{4mm}{28mm}}}\put(35,5){\color{red}{\rule{4mm}{24mm}}}
\put(45,10){\color{red}{\rule{4mm}{20mm}}}\put(25,15){\color{red}{\rule{4mm}{12mm}}}
\put(5,15){\color{red}{\rule{4mm}{4mm}}}\put(10,10){\color{red}{\rule{4mm}{16mm}}}
\put(20,10){\color{red}{\rule{4mm}{16mm}}}\put(40,10){\color{red}{\rule{4mm}{16mm}}}
\put(30,0){\line(1,0){10}} \put(25,5){\line(1,0){25}}
\put(10,10){\line(1,0){40}} \put(0,15){\line(1,0){50}}
\put(0,20){\line(1,0){50}} \put(0,25){\line(1,0){50}}
\put(5,30){\line(1,0){45}} \put(20,35){\line(1,0){30}}
\put(15,30){\line(1,0){10}} \put(25,40){\line(1,0){20}}
\put(0,15){\line(0,1){10}}\put(5,15){\line(0,1){15}}
\put(10,10){\line(0,1){20}}\put(15,10){\line(0,1){20}}\put(20,10){\line(0,1){25}}
\put(25,5){\line(0,1){35}}\put(30,0){\line(0,1){40}}\put(35,0){\line(0,1){40}}\put(40,0){\line(0,1){40}}
\put(45,5){\line(0,1){35}} \put(50,5){\line(0,1){30}}
\put(30,0){\makebox(6,4)[c]{\small \1}}
\put(45,5){\makebox(6,4)[c]{\small \1}}
\put(25,10){\makebox(6,4)[c]{\small \1}}
\put(15,15){\makebox(6,4)[c]{\small\1}}
\put(0,15){\makebox(6,4)[c]{\small \1}}
\put(5,20){\makebox(6,4)[c]{\small\1}}
\put(20,30){\makebox(6,4)[c]{\small\1}}
\put(40,30){\makebox(6,4)[c]{\small \1}}
\put(35,35){\makebox(6,4)[c]{\small \1}}
\end{picture}}
\hspace{2cm}
 {\setlength{\unitlength}{0.8mm}
\begin{picture}(20,50)(0,-5)
\put(0,35){\makebox(6,4)[l]{$c^1=(1,1)$}}
\put(0,30){\makebox(6,4)[l]{$c^2=(0,1,0)$}}
\put(0,25){\makebox(6,4)[l]{$c^3=(0)$}}
\put(0,20){\makebox(6,4)[l]{$c^4=(1,1)$}}
\put(0,15){\makebox(6,4)[l]{$c^5=(0,0,0)$}}
\put(0,10){\makebox(6,4)[l]{$c^6=(1,0)$}}
\put(0,5){\makebox(6,4)[l]{$c^7=(2,0)$}}
\put(0,0){\makebox(6,4)[l]{$c^8=(0,1)$}}
\put(-20,15){\makebox(6,4)[l]{$\longrightarrow$}}
\put(-18,19){\makebox(6,4)[c]{$f$}}
\end{picture}}
\caption{The mapping $f$}
\end{center}
\end{figure}

In order to show that $f$ is bijective, we describe its reverse $g$.
Let ${\bf c}=(c^{(1)},c^{(2)},\ldots,c^{(s)})$ in
$\C_{m_{1}+1}(h_{1})\times\C_{m_{2}+1}(h_{2})\times \cdots\times
\C_{m_{s}+1}(h_{s})$.

Then define the filling $g({\bf c})$ by the following process. We
start with the polyomino (empty filling) $T$.

 (1) Color the columns indexed by the set $A$. We denote by $F_0$ the
result.

(2) Suppose $R_{i_1}\prec R_{i_2} \prec \cdots \prec R_{i_s}$. For
$j$ from 1 to $s$, the (colored) filling $F_{j}$ is obtained from
$F_{j-1}$ by the following process:

\begin{itemize}
\item if $m_{i_j}=0$, then do nothing,
\item else, fill the $i_j\,$-th row of $F_{j-1}$ with $m_{i_j}$ 1's in such a
way that the number of uncolored cells strictly
\begin{itemize}
\item to the left of the first~$1$ is $c^{(i_j)}_1$,
\item between the $u$-th 1 and the $(u+1)$-th 1,
$1\leq u\leq m_{i_j}-1$, is $c^{(i_j)}_{u+1}$,
\item to the right of the last~$1$ is $c^{(i_j)}_{m_{i_j}+1}$.
\end{itemize}
Next, color the cells which are below (resp., above) the new $1$'s
and contained in the $i_j$-th rectangle if $R_{i_j}\in Up(T)$
(resp., $R_{i_j}\in Low(T)$).
\end{itemize}

(3) Set $g({\bf c})=F_s$.

For a better understanding, we give an example. Suppose $T$ is the
moon polyomino given below, $A=\{2\}$ and ${\bf m}=(1,2,1,0,1)$.

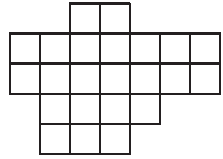
\begin{center}
{\setlength{\unitlength}{0.8mm}
\begin{picture}(35,25)(0,0)
\put(5,0){\line(1,0){15}}\put(5,5){\line(1,0){20}}\put(0,10){\line(1,0){35}}
\put(0,15){\line(1,0){35}}\put(0,20){\line(1,0){35}}\put(10,25){\line(1,0){10}}
\put(0,10){\line(0,1){10}}\put(5,0){\line(0,1){20}}\put(10,0){\line(0,1){25}}
\put(15,0){\line(0,1){25}}\put(20,0){\line(0,1){25}}\put(25,5){\line(0,1){15}}
\put(30,10){\line(0,1){10}}\put(35,10){\line(0,1){10}}
\end{picture}}
\end{center}

Suppose ${\bf c}=(c^{(1)},c^{(2)},c^{(3)},c^{(4)},c^{(5)})$ with
$c^{(1)}=(1,0)$, $c^{(2)}=(1,0,1)$, $c^{(3)}=(0,0,0)$, $c^{(4)}=(0)$
and $c^{(5)}=(0,0)$. The step by step construction of $g({\bf c})$
goes as follows.

\begin{figure}[h!]\label{fig:compo-fillings}
\begin{center}
 {\setlength{\unitlength}{1mm}
\begin{picture}(35,25)(0,0)
\put(-2,3){\makebox(6,4)[c]{$F_0$}}
\put(5,0){\color{red}{\rule{5mm}{20mm}}}
\put(5,0){\line(1,0){15}}\put(5,5){\line(1,0){20}}\put(0,10){\line(1,0){35}}
\put(0,15){\line(1,0){35}}\put(0,20){\line(1,0){35}}\put(10,25){\line(1,0){10}}
\put(0,10){\line(0,1){10}}\put(5,0){\line(0,1){20}}\put(10,0){\line(0,1){25}}
\put(15,0){\line(0,1){25}}\put(20,0){\line(0,1){25}}\put(25,5){\line(0,1){15}}
\put(30,10){\line(0,1){10}}\put(35,10){\line(0,1){10}}
\end{picture}}
\hspace{1cm}
 {\setlength{\unitlength}{1mm}
\begin{picture}(35,25)(0,0)
\put(-2,3){\makebox(6,4)[c]{$F_1$}}
\put(5,0){\color{red}{\rule{5mm}{20mm}}}
\put(15,0){\color{red}{\rule{5mm}{20mm}}}
\put(5,0){\line(1,0){15}}\put(5,5){\line(1,0){20}}\put(0,10){\line(1,0){35}}
\put(0,15){\line(1,0){35}}\put(0,20){\line(1,0){35}}\put(10,25){\line(1,0){10}}
\put(0,10){\line(0,1){10}}\put(5,0){\line(0,1){20}}\put(10,0){\line(0,1){25}}
\put(15,0){\line(0,1){25}}\put(20,0){\line(0,1){25}}\put(25,5){\line(0,1){15}}
\put(30,10){\line(0,1){10}}\put(35,10){\line(0,1){10}}
\put(15,20){\makebox(6,4)[c]{\small \1}}
\put(10,-7){\makebox(6,4)[c]{\small $c^{(i_1)}=(1,0)$}}
\linethickness{0.5mm} \put(10,20){\line(1,0){10}}
\linethickness{0.7mm}
\put(10,0){\line(1,0){10}}\put(10,25){\line(1,0){10}}
\put(10,0){\line(0,1){25}}\put(20,0){\line(0,1){25}}
\end{picture}}
\hspace{1cm}
 {\setlength{\unitlength}{1mm}
\begin{picture}(35,25)(0,0)
\put(-2,3){\makebox(6,4)[c]{$F_2$}}
\put(5,0){\color{red}{\rule{5mm}{20mm}}}
\put(15,0){\color{red}{\rule{5mm}{20mm}}}
\put(10,5){\color{red}{\rule{5mm}{15mm}}}
\put(5,0){\line(1,0){15}}\put(5,5){\line(1,0){20}}\put(0,10){\line(1,0){35}}
\put(0,15){\line(1,0){35}}\put(0,20){\line(1,0){35}}\put(10,25){\line(1,0){10}}
\put(0,10){\line(0,1){10}}\put(5,0){\line(0,1){20}}\put(10,0){\line(0,1){25}}
\put(15,0){\line(0,1){25}}\put(20,0){\line(0,1){25}}\put(25,5){\line(0,1){15}}
\put(30,10){\line(0,1){10}}\put(35,10){\line(0,1){10}}
\put(15,20){\makebox(6,4)[c]{\small \1}}
\put(10,0){\makebox(6,4)[c]{\small \1}}
\put(10,-7){\makebox(6,4)[c]{\small $c^{(i_2)}=(0,0)$}}
\linethickness{0.5mm} \put(5,5){\line(1,0){15}}
\linethickness{0.7mm}
\put(5,0){\line(1,0){15}}\put(5,20){\line(1,0){15}}
\put(5,0){\line(0,1){20}}\put(20,0){\line(0,1){20}}
\end{picture}}
\end{center}
\vspace{1cm}
\begin{center}
 {\setlength{\unitlength}{1mm}
\begin{picture}(35,25)(0,0)
\put(-2,3){\makebox(6,4)[c]{$F_3$}}
\put(5,0){\color{red}{\rule{5mm}{20mm}}}
\put(15,0){\color{red}{\rule{5mm}{20mm}}}
\put(10,5){\color{red}{\rule{5mm}{15mm}}}
\put(5,0){\line(1,0){15}}\put(5,5){\line(1,0){20}}\put(0,10){\line(1,0){35}}
\put(0,15){\line(1,0){35}}\put(0,20){\line(1,0){35}}\put(10,25){\line(1,0){10}}
\put(0,10){\line(0,1){10}}\put(5,0){\line(0,1){20}}\put(10,0){\line(0,1){25}}
\put(15,0){\line(0,1){25}}\put(20,0){\line(0,1){25}}\put(25,5){\line(0,1){15}}
\put(30,10){\line(0,1){10}}\put(35,10){\line(0,1){10}}
\put(15,20){\makebox(6,4)[c]{\small \1}}
\put(10,0){\makebox(6,4)[c]{\small \1}}
\put(10,-7){\makebox(6,4)[c]{\small $c^{(i_3)}=(0)$}}
\linethickness{0.5mm} \put(5,10){\line(1,0){20}}
\linethickness{0.7mm}
\put(5,5){\line(1,0){20}}\put(5,20){\line(1,0){20}}
\put(5,5){\line(0,1){15}}\put(25,5){\line(0,1){15}}
\end{picture}}
\hspace{1cm}
 {\setlength{\unitlength}{1mm}
\begin{picture}(35,25)(0,0)
\put(-2,3){\makebox(6,4)[c]{$F_4$}}
\put(5,0){\color{red}{\rule{5mm}{20mm}}}
\put(15,0){\color{red}{\rule{5mm}{20mm}}}
\put(10,5){\color{red}{\rule{5mm}{15mm}}}
\put(20,10){\color{red}{\rule{5mm}{5mm}}}
\put(25,10){\color{red}{\rule{5mm}{5mm}}}
\put(5,0){\line(1,0){15}}\put(5,5){\line(1,0){20}}\put(0,10){\line(1,0){35}}
\put(0,15){\line(1,0){35}}\put(0,20){\line(1,0){35}}\put(10,25){\line(1,0){10}}
\put(0,10){\line(0,1){10}}\put(5,0){\line(0,1){20}}\put(10,0){\line(0,1){25}}
\put(15,0){\line(0,1){25}}\put(20,0){\line(0,1){25}}\put(25,5){\line(0,1){15}}
\put(30,10){\line(0,1){10}}\put(35,10){\line(0,1){10}}
\put(15,20){\makebox(6,4)[c]{\small \1}}
\put(10,0){\makebox(6,4)[c]{\small \1}}
\put(20,15){\makebox(6,4)[c]{\small \1}}
\put(25,15){\makebox(6,4)[c]{\small \1}}
\put(10,-7){\makebox(6,4)[c]{\small $c^{(i_4)}=(1,0,1)$}}
\linethickness{0.5mm} \put(0,15){\line(1,0){35}}
\linethickness{0.7mm}
\put(0,10){\line(1,0){35}}\put(0,20){\line(1,0){35}}
\put(0,10){\line(0,1){10}}\put(35,10){\line(0,1){10}}
\end{picture}}
\hspace{1cm}
 {\setlength{\unitlength}{1mm}
\begin{picture}(35,25)(0,0)
\put(-2,3){\makebox(6,4)[c]{$F_5$}}
\put(5,0){\color{red}{\rule{5mm}{20mm}}}
\put(15,0){\color{red}{\rule{5mm}{20mm}}}
\put(10,5){\color{red}{\rule{5mm}{15mm}}}
\put(20,10){\color{red}{\rule{5mm}{5mm}}}
\put(25,10){\color{red}{\rule{5mm}{5mm}}}
\put(5,0){\line(1,0){15}}\put(5,5){\line(1,0){20}}\put(0,10){\line(1,0){35}}
\put(0,15){\line(1,0){35}}\put(0,20){\line(1,0){35}}\put(10,25){\line(1,0){10}}
\put(0,10){\line(0,1){10}}\put(5,0){\line(0,1){20}}\put(10,0){\line(0,1){25}}
\put(15,0){\line(0,1){25}}\put(20,0){\line(0,1){25}}\put(25,5){\line(0,1){15}}
\put(30,10){\line(0,1){10}}\put(35,10){\line(0,1){10}}
\put(15,20){\makebox(6,4)[c]{\small \1}}
\put(10,0){\makebox(6,4)[c]{\small \1}}
\put(20,15){\makebox(6,4)[c]{\small \1}}
\put(25,15){\makebox(6,4)[c]{\small \1}}
\put(0,10){\makebox(6,4)[c]{\small \1}}
\put(30,10){\makebox(6,4)[c]{\small \1}}
\put(10,-7){\makebox(6,4)[c]{\small $c^{(i_5)}=(0,0,0)$}}
\linethickness{0.7mm}
\put(0,10){\line(1,0){35}}\put(0,15){\line(1,0){35}}
\put(0,10){\line(0,1){5}}\put(35,10){\line(0,1){5}}
\end{picture}}
\vspace{0.5cm}
 \caption{ The step-by-step construction of $g({\bf
c})$}
\end{center}
\end{figure}
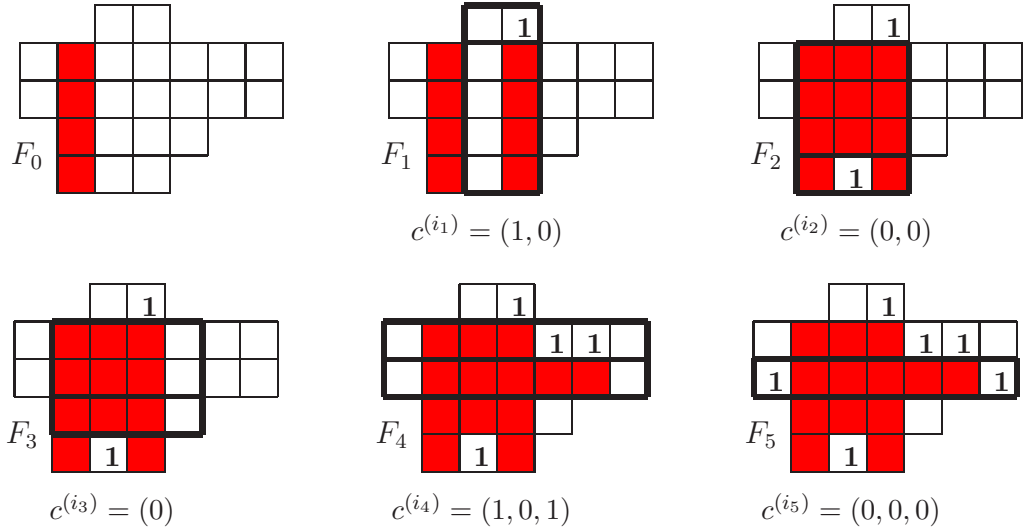

It is not difficult to prove that $g$ is the reverse of $f$, and
thus $f$ is bijective.

Now let ${\bf
c}=(c^{(1)},c^{(2)},\ldots,c^{(s)})\in\C_{m_{1}+1}(h_{1})\times\C_{m_{2}+1}(h_{2})\times
\cdots\times \C_{m_{s}+1}(h_{s})$ and $F=g({\bf c})$ be the
corresponding $01$-filling. Let $i$ be an integer in $[s]$ and $ce$
be the cell of the $i$-th row $R_i$ of $F$ which contains the $j$-th
$1$ of $R_i$. It then follows from the definition of $g$ that
\begin{align*}
\luc(ce;F)&=c_{1}+c_{2}+\cdots+c_{j}\;,\\
\ruc(ce;F)&=c_{j+1}+c_{j+2}+\cdots+c_{m_i+1}=h_i-(c_{1}+c_{2}+\cdots+c_{j}).
\end{align*}
 Applying Proposition~\ref{prop:colored filling} and after some
elementary manipulations of $p,q$-calculus, we get the desired result.\\

\begin{rmk}
Let $c=(c_1,c_2,\ldots,c_k)$ be a composition. Define the
\emph{reverse} $rev(c)$ of $c$ as the composition
$rev(c)=(c_k,c_{k-1},\ldots,c_1)$. Given a sequence of compositions
${\bf c}=(c^{(1)},c^{(2)},\ldots,c^{(s)})$, set $\phi({\bf
c})=(rev(c^{(1)}),rev(c^{(2)}),\ldots,rev(c^{(s)}))$. Then one can
check that the involution $\Phi$ can be factorized as follows:
$$
\Phi=g\circ\phi\circ g^{-1}.
$$
\end{rmk}

\begin{rmk}
It follows from Theorem~\ref{thm:distribution} that the number of
$01$-fillings with at most one $1$ in each column of a given moon
polyomino, such that the number of ascents is $k$ and the number of
descents is $\ell$, does not depend on the order of the rows, given
that the resulting polyomino is again a moon polyomino.
\end{rmk}


\section{Application: symmetry of crossings and nestings
over linked partitions, set partitions and matchings}

  Let $G$ be a graph on $[n]$. The multiset of lefthand (resp.,
righthand) endpoints of the arcs of $G$ will be denoted by $\le(G)$
(resp., $\ri(G)$). For example, if $G$ is the graph given in
Figure~1, we have $\le(G) = \{1,2,2,3, 5,6,6,9\}$ and $\ri(G) =
\{3,4, 6, 7, 9,9,10,11\}$. Suppose $F$ is the $01$-filling of the
triangular shape $\Delta_n$ associated to $G$. Then it is easy to
see that the number of $1$'s in the column (resp., row) labelled~$i$
is equal to the multiciplicity of $i$ in $\le(G)$ (resp., $\ri(G)$).
See Figure~1 for an example.
As explained in the introduction, results on $2$-crossings and
$2$-nestings in simple graphs translate into results on ascents and
descents in fillings of  $\Delta_n$ and vice-versa. It is then easy
to derive the following result from Theorem~\ref{thm:refinement
main1}.

 If $S$ and $T$ are two
multisubsets of $[n]$, denote by $\G_n(S,T)$ the set of graphs $G$
on $[n]$ satisfying $\le(G)=S$ and $\ri(G)=T$.
\begin{thm}\label{thm:symgraphes}
For any pair $(S,T)$ of multisubsets of $[n]$ such that
\begin{enumerate}
\item either all elements of $S$ have multiplicity $1$,
\item either all elements of $T$ have multiplicity $1$,
\end{enumerate}
 the joint statistic $(\cro_2,\ne_2)$ is
symmetrically distributed over $\G_n(S,T)$.
\end{thm}

 It is now easy to recover $\eqref{eq:sympartperm}$.

\subsection{Linked partitions}
 Let $E$ and $F$ be two finite subsets of integers. We say that $E$ and
$F$ are \emph{nearly disjoint} if for every $i\in E\cap F$, one of
the following holds:
\begin{itemize}
  \item[(a)] $i = \min(E)$, $|E| > 1$ and $i\neq \min(F)$, or
  \item[(b)] $i = \min(F)$, $|F| > 1$ and $i\neq \min(E)$.
\end{itemize}
 A \emph{linked partition} (see \cite{Chen1}) of
 $[n]$ is a collection of nonempty and
pairwise nearly disjoint subsets  whose union is $[n]$.
 The set of all linked partitions
of $[n]$ will be denoted by $\LP_n$. The linear representation
$G_{\pi}$ of a linked partition $\pi\in\LP_n$ is the graph on $[n]$
where $i$ and $j$ are connected by an arc if and only if $j$ lies in
a block B with $i= \Min(B)$. An illustration is given in Figure~.
Clearly, this establishes a bijection between linked set partitions
and simple graphs $G$ such that all elements of $\ri(G)$ have
multiplicity one.

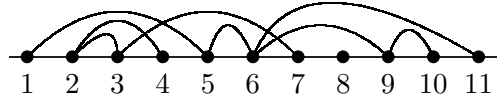
\begin{figure}[h!]
\begin{center}
{\setlength{\unitlength}{1.2mm}
\begin{picture}(50,13)(0,-3)
\put(-2,0){\line(1,0){54}}
\put(0,0){\circle*{1,3}}\put(0,0){\makebox(0,-6)[c]{\small $1$}}
\put(5,0){\circle*{1,3}}\put(5,0){\makebox(0,-6)[c]{\small $2$}}
\put(10,0){\circle*{1,3}}\put(10,0){\makebox(0,-6)[c]{\small $3$}}
\put(15,0){\circle*{1,3}}\put(15,0){\makebox(0,-6)[c]{\small $4$}}
\put(20,0){\circle*{1,3}}\put(20,0){\makebox(0,-6)[c]{\small $5$}}
\put(25,0){\circle*{1,3}}\put(25,0){\makebox(0,-6)[c]{\small $6$}}
\put(30,0){\circle*{1,3}}\put(30,0){\makebox(0,-6)[c]{\small $7$}}
\put(35,0){\circle*{1,3}}\put(35,0){\makebox(0,-6)[c]{\small $8$}}
\put(40,0){\circle*{1,3}}\put(40,0){\makebox(0,-6)[c]{\small $9$}}
\put(45,0){\circle*{1,3}}\put(45,0){\makebox(0,-6)[c]{\small $10$}}
\put(50,0){\circle*{1,3}}\put(50,0){\makebox(0,-6)[c]{\small $11$}}
\qbezier(0,0)(10,10)(20,0)
\qbezier(5,0)(10,5)(10,0)\qbezier(20,0)(22,7)(25,0)
\qbezier(10,0)(20,10)(30,0)\qbezier(25,0)(30,12)(50,0)
\qbezier(25,0)(32.5,7)(40,0)\qbezier(40,0)(42,6)(45,0)
\qbezier(5,0)(10,8)(15,0)
\end{picture}}
\end{center}
\caption{Linear representation of
$\pi=\{1,5\}\{2,3,4\}\{3,7\}\{5,6\}\{6,9,11\}\{8\}\{9,10\}$}
\end{figure}

For $\O,\C \subseteq[n]$ two multisubsets of $[n]$, denote by
$\LP_n(\O,\C)$ the set $\{\pi\in\LP_n\,:\, \le(G_{\pi})=
\O\,,\,\ri(G_{\pi}) = \C\}$. We then derive from
Theorem~\ref{thm:symgraphes} the following result due to Chen et al
\cite{Chen1}.

\begin{cor}\label{cor:chen}
For any integer $n$ and for any pair $(\O,\C)$ of multisubsets  of
$[n]$, we have
\begin{equation*}\label{eq:sym-linked}
\sum_{\pi\in\,\LP_n(\O,\C)}p^{\cro_2(G_{\pi})}q^{\ne_2(G_{\pi})}
=\sum_{\pi\in\,\LP_n(\O,\C)}p^{\ne_2(G_{\pi})}q^{\cro_2(G_{\pi})}.
\end{equation*}
Summing over all $(\O,\C)$, we get
\begin{equation*}\label{eq:sym-linked}
\sum_{\pi\in\,\LP_n}p^{\cro_2(G_{\pi})}q^{\ne_2(G_{\pi})}
=\sum_{\pi\in\,\LP_n}p^{\ne_2(G_{\pi})}q^{\cro_2(G_{\pi})}.
\end{equation*}
\end{cor}

\subsection{Matchings and set partitions}
Recall that a (\emph{set}) \emph{partition} of $[n]$ is a collection
of nonempty pairwise-disjoint sets, called blocks, whose union is
$[n]$. A (\emph{complete}) \emph{matching} is just a partition whose
each block contains exactly two elements. In fact, set partitions
and matchings are just particular linked partitions. The set of all
set partitions and matchings of $[n]$ will be denoted  respectively
by $\P_n$ and $\M_n$. Set partitions (and thus matchings) have a
natural graphical representation, called \emph{standard
repesentation}. To each set partition $\pi$ of $[n]$, one associates
the graph $St_{\pi}$ on $[n]$ whose edge set consists of arcs
joining the elements of each block in numerical order. An
illustration is given in Figure~8. This establishes a bijection
between set partitions and simple graphs $G$ such that all elements
of $\ri(G)$ and $\le(G)$ have multiplicity one. This also
establishes a bijection matchings and simple graphs $G$ such that
$\ri(G)\cap\le(G)=\emptyset$, $\ri(G)\cup\le(G)=[n]$ and all
elements of $\ri(G)$ and $\le(G)$ have multiplicity one. Note that
for $M\in\M_n$, we have $St_M=G_M$ but if $\pi\in\P_n$ contains a
block of length at least $3$ then $St_{\pi}\neq G_{\pi}$.

\begin{figure}[h!]
\begin{center}
{\setlength{\unitlength}{1mm}
\begin{picture}(50,10)(0,0)
\put(-2,0){\line(1,0){49}}
\put(0,0){\circle*{1,3}}\put(0,0){\makebox(0,-6)[c]{\small 1}}
\put(5,0){\circle*{1,3}}\put(5,0){\makebox(0,-6)[c]{\small 2}}
\put(10,0){\circle*{1,3}}\put(10,0){\makebox(0,-6)[c]{\small 3}}
\put(15,0){\circle*{1,3}}\put(15,0){\makebox(0,-6)[c]{\small 4}}
\put(20,0){\circle*{1,3}}\put(20,0){\makebox(0,-6)[c]{\small 5}}
\put(25,0){\circle*{1,3}}\put(25,0){\makebox(0,-6)[c]{\small 6}}
\put(30,0){\circle*{1,3}}\put(30,0){\makebox(0,-6)[c]{\small 7}}
\put(35,0){\circle*{1,3}}\put(35,0){\makebox(0,-6)[c]{\small 8}}
\put(40,0){\circle*{1,3}}\put(40,0){\makebox(0,-6)[c]{\small 9}}
\put(45,0){\circle*{1,3}}\put(45,0){\makebox(0,-6)[c]{\small 10}}
\qbezier(0,0)(10,8)(20,0) \qbezier(25,0)(32,6)(40,0)
\qbezier(5,0)(20,12)(35,0)\qbezier(15,0)(30,10)(45,0)
\qbezier(10,0)(20,8)(30,0)
\end{picture}
}\hspace{2cm} {\setlength{\unitlength}{1mm}
\begin{picture}(50,10)(0,0)
\put(-2,0){\line(1,0){54}}
\put(0,0){\circle*{1,3}}\put(0,0){\makebox(0,-6)[c]{\small 1}}
\put(5,0){\circle*{1,3}}\put(5,0){\makebox(0,-6)[c]{\small 2}}
\put(10,0){\circle*{1,3}}\put(10,0){\makebox(0,-6)[c]{\small 3}}
\put(15,0){\circle*{1,3}}\put(15,0){\makebox(0,-6)[c]{\small 4}}
\put(20,0){\circle*{1,3}}\put(20,0){\makebox(0,-6)[c]{\small 5}}
\put(25,0){\circle*{1,3}}\put(25,0){\makebox(0,-6)[c]{\small 6}}
\put(30,0){\circle*{1,3}}\put(30,0){\makebox(0,-6)[c]{\small 7}}
\put(35,0){\circle*{1,3}}\put(35,0){\makebox(0,-6)[c]{\small 8}}
\put(40,0){\circle*{1,3}}\put(40,0){\makebox(0,-6)[c]{\small 9}}
\put(45,0){\circle*{1,3}}\put(45,0){\makebox(0,-6)[c]{\small 10}}
\put(50,0){\circle*{1,3}}\put(50,0){\makebox(0,-6)[c]{\small 11}}
\qbezier(0,0)(20,12)(40,0) \qbezier(40,0)(42,5)(45,0)
\qbezier(5,0)(7,5)(10,0)\qbezier(20,0)(22,5)(25,0)
\qbezier(10,0)(20,8)(30,0)\qbezier(25,0)(37,10)(50,0)
\end{picture}
}
\end{center}
\caption{\small{Standard representations of
$M=\{1,5\}\{2,8\}\{3,7\}\{4,10\}\{6,9\}$ and
$\pi=\{1,9,10\}\{2,3,7\}\{4\}\{5,6,11\}\{8\}$}}
\end{figure}
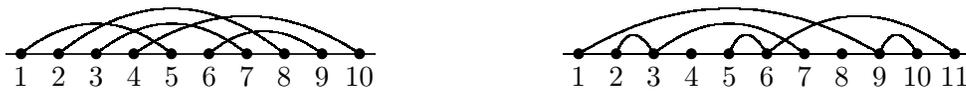

For $\O,\C \subseteq[n]$ two subsets of $[n]$, let $\P_n(\O,\C)$ be
the sets of set partitions $\{\pi\in\P_n\,:\, \le(G_{\pi})=
\O\,,\,\ri(G_{\pi}) = \C\}$ and $\M_n(\O,\C)$ the set of matchings
$\{M\in\M_n\,:\, \le(G_{M})= \O\,,\,\ri(G_{M}) = \C\}$. Then we
derive immediately from Theorem~\ref{thm:symgraphes} (or
Corollary~\ref{cor:chen}) the following results whose first of them
is due to Kasraoui and Zeng \cite{KaZe} and second to Klazar and Noy
\cite{KlNo}.

\begin{cor}
For any integer $n$ and pair $(\O,\C)$ of subsets of $[n]$, we have
\begin{equation*}\label{eq:sym-part}
\sum_{\pi\in\,\P_n(\O,\C)}p^{\cro_2(St_{\pi})}q^{\ne_2(St_{\pi})}
=\sum_{\pi\in\,\P_n(\O,\C)}p^{\ne_2(St_{\pi})}q^{\cro_2(St_{\pi})}.
\end{equation*}
Summing over all $(\O,\C)$, we get
\begin{equation*}\label{eq:sym-partgen}
\sum_{\pi\in\,\P_n}p^{\cro_2(St_{\pi})}q^{\ne_2(St_{\pi})}
=\sum_{\pi\in\,\P_n}p^{\ne_2(St_{\pi})}q^{\cro_2(St_{\pi})}.
\end{equation*}
\end{cor}

\begin{cor}
For any integer $n$ and pair $(\O,\C)$ of subsets of $[n]$, we have
\begin{equation*}\label{eq:sym-mat}
\sum_{M\in\,\M_n(\O,\C)}p^{\cro_2(St_{M})}q^{\ne_2(St_{M})}
=\sum_{M\in\,\M_n(\O,\C)}p^{\ne_2(St_{M})}q^{\cro_2(St_{M})}.
\end{equation*}
Summing over all $(\O,\C)$, we get
\begin{equation*}\label{eq:sym-matgen}
\sum_{M\in\,\M_n}p^{\cro_2(St_{M})}q^{\ne_2(St_{M})}
=\sum_{M\in\,\M_n}p^{\ne_2(St_{M})}q^{\cro_2(St_{M})}.
\end{equation*}
\end{cor}

\subsection{Distribution of $(\cro_2,\ne_2)$ over $\LP_n(\O,\C)$, $\P_n(\O,\C)$
and $\M_n(\O,\C)$} We just want to point briefly that the
distribution of $(\cro_2,\ne_2)$ over $\LP_n(\O,\C)$, $\P_n(\O,\C)$
and $\M_n(\O,\C)$ can be easily derived from
Theorem~\ref{thm:distribution} and the correspondence between simple
graphs $[n]$ and $01$-fillings of $\Delta_n$.

Let $(\O,\C)$ be a pair of multisubsets of $[n]$ and denote by $m_i$
the multiplicity of $i\in\O$. Also for any $i\in\O$ set
$h_i=|\{j\in\C\,|\,j>i\}|-|\{j\in\O\,|\,j>i\}|$. Then it can be
deduced from Theorem~\ref{thm:distribution} that
\begin{equation}
\sum_{\pi\in\,\LP_n(\O,\C)}p^{\cro_2(G_{\pi})}q^{\ne_2(G_{\pi})}
=\prod_{i\in\O}{h_i\brack m_i}_{p,q}.
\end{equation}
Note that the above identity is equivalent to a result of Chen et al
\cite[Theorem~3.5]{Chen1}. In set partitions, $m_i\leq 1$ for any
$i\in\O$. We then get the following result which is implicit
in~\cite[Section~4]{KaZe}
\begin{equation}
\sum_{\pi\in\,\P_n(\O,\C)}p^{\cro_2(G_{\pi})}q^{\ne_2(G_{\pi})}
=\prod_{i\in\O}[h_i]_{p,q}.
\end{equation}

\section{Concluding remarks}

\subsection{}
 It is natural, in view of Theorem~\ref{thm: main}, to ask  if
the symmetry of the joint distribution of the bi-statistic
$(\nec_2,\sec_2)$ extend for arbitrary $01$-fillings of moon
polyominoes, i.e., no restrictions on the number of $1$'s in columns
and rows. The answer is no by means of the following result.

\begin{prop}
For any $n\geq5$ the numbers of arbitrary $01$-fillings
of~$\Delta_n$
\begin{itemize}
\item with exactly ${n\choose 4}$ descents is equal to $2^n$,
\item with exactly ${n\choose 4}$ ascents is equal to $16$.
\end{itemize}
This implies that   the joint distribution of $(\nec_2,\sec_2)$ over
all arbitrary $01$-fillings of $\Delta_n$ is not symmetric for any
$n\geq 5$.
\end{prop}

\pf We give the proof for $n=5,6$ since the reasoning can be
generalized for arbitrary $n$. Suppose $n=5$. Then one can check
that the arbitrary $01$-fillings of $\Delta_5$ with exactly $5$
descents and those with exactly $5$ ascents have the following
"form"
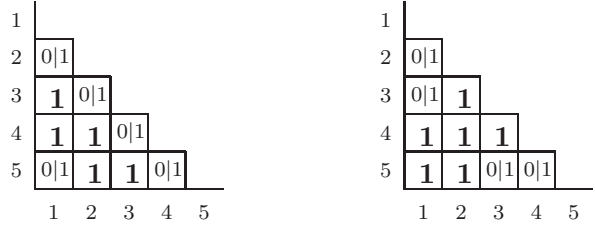
\begin{center}
{\setlength{\unitlength}{1mm}
\begin{picture}(25,30)(0,-5)
\put(0,0){\line(1,0){25}}\put(0,5){\line(1,0){20}}\put(0,10){\line(1,0){15}}
\put(0,15){\line(1,0){10}}\put(0,20){\line(1,0){5}}
\put(0,0){\line(0,1){25}}\put(5,0){\line(0,1){20}}\put(10,0){\line(0,1){15}}
\put(15,0){\line(0,1){10}}\put(20,0){\line(0,1){5}}
\put(0,5){\makebox(6,4)[c]{\small \1}}
\put(0,10){\makebox(6,4)[c]{\small \1}}
\put(5,0){\makebox(6,4)[c]{\small \1}}
\put(5,5){\makebox(6,4)[c]{\small \1}}
\put(10,0){\makebox(6,4)[c]{\small \1}}
\put(0,0){\makebox(5.5,5)[c]{\tiny $0|1$}}
\put(0,15){\makebox(5.5,5)[c]{\tiny $0|1$}}
\put(5,10){\makebox(5.5,5)[c]{\tiny $0|1$}}
\put(10,5){\makebox(5.5,5)[c]{\tiny $0|1$}}
\put(15,0){\makebox(5.5,5)[c]{\tiny $0|1$}}
\put(0,0){\makebox(5,-6)[c]{\tiny $1$}}
\put(5,0){\makebox(5,-6)[c]{\tiny $2$}}
\put(10,0){\makebox(5,-6)[c]{\tiny $3$}}
\put(15,0){\makebox(5,-6)[c]{\tiny $4$}}
\put(20,0){\makebox(5,-6)[c]{\tiny $5$}}
\put(0,0){\makebox(-5,5)[c]{\tiny $5$}}
\put(0,5){\makebox(-5,5)[c]{\tiny $4$}}
\put(0,10){\makebox(-5,5)[c]{\tiny $3$}}
\put(0,15){\makebox(-5,5)[c]{\tiny $2$}}
\put(0,20){\makebox(-5,5)[c]{\tiny $1$}}
\end{picture}}
\hspace{2cm} {\setlength{\unitlength}{1mm}
\begin{picture}(25,30)(0,-5)
\put(0,0){\line(1,0){25}}\put(0,5){\line(1,0){20}}\put(0,10){\line(1,0){15}}
\put(0,15){\line(1,0){10}}\put(0,20){\line(1,0){5}}
\put(0,0){\line(0,1){25}}\put(5,0){\line(0,1){20}}\put(10,0){\line(0,1){15}}
\put(15,0){\line(0,1){10}}\put(20,0){\line(0,1){5}}
\put(0,0){\makebox(6,4)[c]{\small \1}}
\put(0,5){\makebox(6,4)[c]{\small \1}}
\put(5,0){\makebox(6,4)[c]{\small \1}}
\put(5,5){\makebox(6,4)[c]{\small \1}}
\put(5,10){\makebox(6,4)[c]{\small \1}}
\put(10,5){\makebox(6,4)[c]{\small \1}}
\put(0,10){\makebox(5.5,5)[c]{\tiny $0|1$}}
\put(0,15){\makebox(5.5,5)[c]{\tiny $0|1$}}
\put(10,0){\makebox(5.5,5)[c]{\tiny $0|1$}}
\put(15,0){\makebox(5.5,5)[c]{\tiny $0|1$}}
\put(0,0){\makebox(5,-6)[c]{\tiny $1$}}
\put(5,0){\makebox(5,-6)[c]{\tiny $2$}}
\put(10,0){\makebox(5,-6)[c]{\tiny $3$}}
\put(15,0){\makebox(5,-6)[c]{\tiny $4$}}
\put(20,0){\makebox(5,-6)[c]{\tiny $5$}}
\put(0,0){\makebox(-5,5)[c]{\tiny $5$}}
\put(0,5){\makebox(-5,5)[c]{\tiny $4$}}
\put(0,10){\makebox(-5,5)[c]{\tiny $3$}}
\put(0,15){\makebox(-5,5)[c]{\tiny $2$}}
\put(0,20){\makebox(-5,5)[c]{\tiny $1$}}
\end{picture}}
\end{center}

from which it is easy to obtain the result. Similarly, for $n=6$,
the arbitrary $01$-fillings of $\Delta_6$ with exactly $15$ descents
and those with exactly $15$ ascents have the following "form".

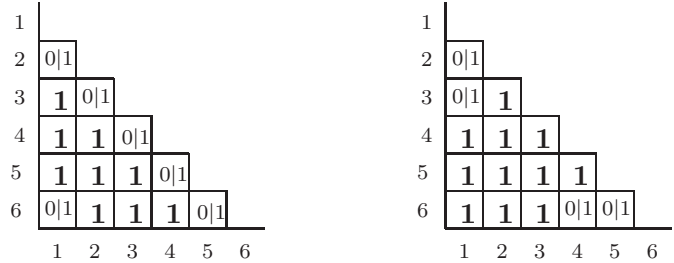
\begin{center}
{\setlength{\unitlength}{1mm}
\begin{picture}(30,35)(0,-5)
\put(0,0){\line(1,0){30}}\put(0,5){\line(1,0){25}}\put(0,10){\line(1,0){20}}
\put(0,15){\line(1,0){15}}\put(0,20){\line(1,0){10}}\put(0,25){\line(1,0){5}}
\put(0,0){\line(0,1){30}}\put(5,0){\line(0,1){25}}\put(10,0){\line(0,1){20}}
\put(15,0){\line(0,1){15}}\put(20,0){\line(0,1){10}}\put(25,0){\line(0,1){5}}
\put(0,0){\makebox(5.5,5)[c]{\tiny $0|1$}}
\put(0,5){\makebox(6,4)[c]{\small \1}}
\put(0,10){\makebox(6,4)[c]{\small \1}}
\put(0,15){\makebox(6,4)[c]{\small \1}}
\put(0,20){\makebox(5.5,5)[c]{\tiny $0|1$}}
\put(5,0){\makebox(6,4)[c]{\small \1}}
\put(5,5){\makebox(6,4)[c]{\small \1}}
\put(5,10){\makebox(6,4)[c]{\small \1}}
\put(5,15){\makebox(5.5,5)[c]{\tiny $0|1$}}
\put(10,0){\makebox(6,4)[c]{\small \1}}
\put(10,5){\makebox(6,4)[c]{\small \1}}
\put(10,10){\makebox(6,4)[c]{\tiny $0|1$}}
\put(15,0){\makebox(6,4)[c]{\small \1}}
\put(15,5){\makebox(6,4)[c]{\tiny $0|1$}}
\put(20,0){\makebox(6,4)[c]{\tiny $0|1$}}
\put(0,0){\makebox(5,-6)[c]{\tiny $1$}}
\put(5,0){\makebox(5,-6)[c]{\tiny $2$}}
\put(10,0){\makebox(5,-6)[c]{\tiny $3$}}
\put(15,0){\makebox(5,-6)[c]{\tiny $4$}}
\put(20,0){\makebox(5,-6)[c]{\tiny $5$}}
\put(25,0){\makebox(5,-6)[c]{\tiny $6$}}
\put(0,0){\makebox(-6,5)[c]{\tiny $6$}}
\put(0,5){\makebox(-6,5)[c]{\tiny $5$}}
\put(0,10){\makebox(-5,5)[c]{\tiny $4$}}
\put(0,15){\makebox(-5,5)[c]{\tiny $3$}}
\put(0,20){\makebox(-5,5)[c]{\tiny $2$}}
\put(0,25){\makebox(-5,5)[c]{\tiny $1$}}
\end{picture}}
\hspace{2cm} {\setlength{\unitlength}{1mm}
\begin{picture}(30,35)(0,-5)
\put(0,0){\line(1,0){30}}\put(0,5){\line(1,0){25}}\put(0,10){\line(1,0){20}}
\put(0,15){\line(1,0){15}}\put(0,20){\line(1,0){10}}\put(0,25){\line(1,0){5}}
\put(0,0){\line(0,1){30}}\put(5,0){\line(0,1){25}}\put(10,0){\line(0,1){20}}
\put(15,0){\line(0,1){15}}\put(20,0){\line(0,1){10}}\put(25,0){\line(0,1){5}}
\put(0,0){\makebox(6,4)[c]{\small \1}}
\put(0,5){\makebox(6,4)[c]{\small \1}}
\put(0,10){\makebox(6,4)[c]{\small \1}}
\put(0,15){\makebox(5.5,5)[c]{\tiny $0|1$}}
\put(0,20){\makebox(5.5,5)[c]{\tiny $0|1$}}
\put(5,0){\makebox(6,4)[c]{\small \1}}
\put(5,5){\makebox(6,4)[c]{\small \1}}
\put(5,10){\makebox(6,4)[c]{\small \1}}
\put(5,15){\makebox(6,4)[c]{\small \1}}
\put(10,0){\makebox(6,4)[c]{\small \1}}
\put(10,5){\makebox(6,4)[c]{\small \1}}
\put(10,10){\makebox(6,4)[c]{\small \1}}
\put(15,5){\makebox(6,4)[c]{\small \1}}
\put(15,0){\makebox(5.5,5)[c]{\tiny $0|1$}}
\put(20,0){\makebox(5.5,5)[c]{\tiny $0|1$}}
\put(0,0){\makebox(5,-6)[c]{\tiny $1$}}
\put(5,0){\makebox(5,-6)[c]{\tiny $2$}}
\put(10,0){\makebox(5,-6)[c]{\tiny $3$}}
\put(15,0){\makebox(5,-6)[c]{\tiny $4$}}
\put(20,0){\makebox(5,-6)[c]{\tiny $5$}}
\put(25,0){\makebox(5,-6)[c]{\tiny $6$}}
\put(0,0){\makebox(-6,5)[c]{\tiny $6$}}
\put(0,5){\makebox(-6,5)[c]{\tiny $5$}}
\put(0,10){\makebox(-5,5)[c]{\tiny $4$}}
\put(0,15){\makebox(-5,5)[c]{\tiny $3$}}
\put(0,20){\makebox(-5,5)[c]{\tiny $2$}}
\put(0,25){\makebox(-5,5)[c]{\tiny $1$}}
\end{picture}}
\end{center}

\qed

 We then propose the following problem.
\begin{problem}
Characterize the moon polyomino $T$ for which the joint distribution
of $(\nec_2,\sec_2)$ over arbitrary $01$-fillings of $T$ is
symmetric?
\end{problem}
Note that $T=\Delta_4$ satisfies the above condition.

\subsection{}

One can also ask if Theorem~\ref{thm:refinement main1} and
Theorem~\ref{thm: main} can be extended to arbitrary larger classes
of polyominoes. We note that the condition of intersection free is
necessary. Indeed, the polyomino $T$ represented below is convex but
not intersection free,
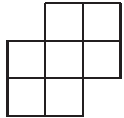
\begin{center}
{\setlength{\unitlength}{1mm}
\begin{picture}(15,20)(0,0)
\put(0,0){\line(1,0){10}}\put(0,5){\line(1,0){15}}\put(0,10){\line(1,0){15}}
\put(5,15){\line(1,0){10}}
\put(0,0){\line(0,1){10}}\put(5,0){\line(0,1){15}}\put(10,0){\line(0,1){15}}
\put(15,5){\line(0,1){10}}
\end{picture}}
\end{center}
and
$$\sum_{F\in\N(T,(1,1,1),\emptyset)}p^{\nec_2(F)}q^{\sec_2(F)}=\sum_{F\in\N(T,(1,1,1))}p^{\nec_2(F)}q^{\sec_2(F)}=p^2+2q$$
is not symmetric. Also note that one can check that $\sum_{\bf
m}\sum_{F\in\N(T,{\bf m})}p^{\nec_2(F)}q^{\sec_2(F)}$ is not
symmetric.


\renewcommand{\baselinestretch}{1}

\end{document}